\documentclass[12pt,reqno]{amsart}
\usepackage[cp1251]{inputenc}
\usepackage[T2A]{fontenc}
\usepackage[english,russian,ukrainian]{babel}
\usepackage{amsmath,amsfonts,amssymb}
\usepackage{geometry}

\numberwithin{equation}{section}

\newtheorem{theorem}{Теорема}[section]
\newtheorem{corollary}{Наслідок}[section]

\theoremstyle{remark}
\newtheorem{remark}{Зауваження}[section]

\newtheorem{example}{Приклад}[section]

\begin{document}

\begin{flushleft}

\textbf{V. A. Mikhailets, A. A. Murach, I. S. Chepurukhina}\\

\small(Institute of Mathematics of the National Academy of Sciences of Ukraine, Kyiv)

\medskip

\large\textbf{ELLIPTIC OPERATORS AND BOUNDARY PROBLEMS\\ IN SPACES OF GENERALIZED SMOTHNESS}

\normalsize

\bigskip

\textbf{В. А. Михайлець, О. О. Мурач, І. С. Чепурухіна}\\

\small{(Інститут математики Національної академії наук України, Київ)}

\medskip

\large\textbf{ЕЛІПТИЧНІ ОПЕРАТОРИ І КРАЙОВІ ЗАДАЧІ\\ У ПРОСТОРАХ УЗАГАЛЬНЕНОЇ ГЛАДКОСТІ}

\end{flushleft}

%%2020 Mathematics Subject Classification: 46E35, 35J40, 35J48.

%% Keywords: space of generalized smoothness, H\"ormander space, elliptic operator, elliptic problem.

\bigskip

\noindent\small The paper contains a survey of the results obtained during the last ten years in the theory of elliptic boundary problems in H\"ormander function spaces, developed by the authors, and other related results of modern analysis. The basics of this theory and some of its applications are systematically expounded in the monograph  \textit{H\"ormander Spaces, Interpolation, and Elliptic Problems} (De Gruyter, Berlin/Boston, 2014) by the first two authors of the survey.

\bigskip

\normalsize

\section{Вступ.}\label{survey-sect1}

У теорії рівнянь з частинними похідними центральне місце посідають питання існування, єдиності та регулярності їх розв'язків. На відміну від звичайних диференціальних рівнянь з регулярними коефіцієнтами ці питання є досить складними. Так, відомі приклади лінійних диференціальних рівнянь з гладкими коефіцієнтами та правими частинами в $\mathbb{R}^{3}$, які не мають розв'язків на жодній відкритій підмножині простору $\mathbb{R}^{3}$ навіть у класі розподілів \cite[п.~6.1]{Hermander63}. Деякі однорідні диференціальні рівняння, зокрема еліптичного типу, з гладкими (але не аналітичними) коефіцієнтами в $\mathbb{R}^{n}$ мають нетривіальні розв'язки з компактним носієм \cite[теорема~13.6.15]{Hermander05v2}. Тому будь-яка крайова задача для такого рівняння в евклідовій області, ширшій за носій, не є однозначно розв'язною. Питання щодо регулярності розв'язків є також досить складним. Наприклад, вже для оператора Лапласа відомо \cite[п.~4.5]{GilbargTrudinger98}, що з умови $\Delta u\in C(\Omega)$ для розподілу $u$ не випливає, що він є двічі диференційовною класичною функцією в евклідовій області $\Omega$. Отже, максимальна регулярність розв'язку рівняння $\Delta u=f$ не досягається у шкалі просторів $C^{k}(\Omega)$, де ціле $k\geq2$, тобто $f\in C^{k-2}(\Omega)\not\Rightarrow u\in C^{k}(\Omega)$. Аналогічна ситуація має місце для довільного лінійного еліптичного рівняння з коефіцієнтами класу $C^{\infty}(\overline{\Omega})$ \cite[п~5.7.3]{Triebel95}.

Найбільш повно зазначені питання досліджені для лінійних еліптичних диференціальних рівнянь і крайових задач. Це було зроблено в 50--60-х роках XX століття у роботах С.~Агмона, А.~Дугліса, Л.~Ніренберга \cite{AgmonDouglisNirenberg59, AgmonDouglisNirenberg64, DouglisNirenberg55}, Ю.~М.~Березанського \cite{Berezansky68}, Ф.~Браудера \cite{Browder56, Browder59}, Л.~Р.~Волевича \cite{Volevich65}, Ж.-Л.~Ліонса і Е.~Мадженеса \cite{LionsMagenes62V, LionsMagenes63VI, LionsMagenes72}, Я.~А.~Ройтберга \cite{Roitberg96, Roitberg65},
В.~О.~Солоннікова \cite{Solonnikov64, Solonnikov66},  Л.~Хермандера \cite{Hermander63}, М.~Шехтера \cite{Schechter60, Schechter63} та інших. При цьому еліптичні рівняння були досліджені у класичних шкалах просторів Гельдера (дробових порядків) і просторів Соболєва (як додатних, так і від'ємних дійсних порядків).

Фундаментальний результат теорії еліптичних рівнянь полягає у тому, що вони породжують обмежені фредгольмові оператори (тобто оператори зі скінченним індексом) на відповідних парах цих просторів. Наприклад, якщо на замкненому гладкому многовиді $\Gamma$ задано лінійне еліптичне диференціальне рівняння $Au=f$ порядку $m$, то диференціальний оператор $A$ є обмеженим і фредгольмовим на парі гільбертових просторів Соболєва $H^{s+m}(\Gamma)$ і $H^{s}(\Gamma)$ відповідно порядків $s+m$ і $s\in\mathbb{R}$. Звідси випливає, що максимальна гладкість розв'язку такого рівняння досягається у соболєвській шкалі, тобто $f\in H^{s}(\Gamma)\Rightarrow u\in H^{s+m}(\Gamma)$. Якщо многовид має межу, то фредгольмів оператор породжується неоднорідною еліптичною крайовою задачею для вказаного рівняння, наприклад, задачею з крайовими умовами Діріхле.
Згодом у працях Ґ.~Трібеля \cite{Triebel83}, Й.~Франке і Т.~Рунста \cite{FrankeRunst95, RunstSickel96}, Й.~Йонсона \cite{Johnsen96} та одного з авторів \cite{Murach94UMJ12, Murach94Dop12} було показано, що еліптичні рівняння і крайові задачі породжують фредгольмові оператори на більш тонких шкалах функціональних просторів Нікольського--Бєсова, Зігмунда і Лізоркіна--Трібеля. Ці результати знайшли численні застосування у теорії диференціальних рівнянь, математичній фізиці, спектральній теорії диференціальних операторів (див. монографії Ю.~М.~Березанського \cite{Berezansky68}, О.~О.~Ладиженської і Н.~М.~Уральцевої \cite{LadyszenskayaUraltseva64}, Ж.-Л.~Ліонса і Е.~Мадженеса \cite{LionsMagenes72}, Я.~А.~Ройтберга \cite{Roitberg96, Roitberg99}, І.~В.~Скрипника \cite{Skrypnik94}, Ґ.~Трібеля \cite{Triebel95}, огляди М.~С.~Аграновича \cite{Agranovich94, Agranovich97} і наведену там бібліографію).

З точки зору застосувань, особливо до спектральної теорії, найбільш важливим є випадок гільбертових просторів. Втім, єдиною шкалою гільбертових просторів, у якій були систематично досліджені властивості еліптичних операторів, тривалий час залишалася шкала просторів Соболєва. Як з'ясувалося вона є занадто грубою для багатьох важливих задач аналізу \cite{Haroske07, Stepanets05, Triebel01, Triebel10}, теорії диференціальних рівнянь \cite{Hermander63, Hermander05v2, NicolaRodino10, Paneah00}, теорії випадкових процесів \cite{Jacob010205}, де виникають еліптичні оператори.

Цього недоліку позбавлені функціональні простори узагальненої гладкості, зокрема введені Л.~Хермандером в \cite[п.~2.2]{Hermander63}. Вони параметризовані не числовим набором, а досить загальною ваговою функцією $k=k(\xi)$ частотних змінних $\xi=(\xi_1,\ldots,\xi_n)$, дуальних до просторових змінних відносно перетворення Фур'є $\mathcal{F}$. Зокрема, Л.~Хермандер дослідив гільбертів простір $\mathcal{B}_{2,k}(\mathbb{R}^{n})$, який складається
з усіх повільно зростаючих розподілів $w$ на $\mathbb{R}^{n}$ таких, що
$k\mathcal{F}w\in L_{2}(\mathbb{R}^{n})$, і наділений нормою $\|k\mathcal{F}w\|_{L_{2}(\mathbb{R}^{n})}$ розподілу $w$. Якщо $k(\xi)\equiv(1+|\xi|^{2})^{s/2}$ для деякого $s\in\mathbb{R}$, то отримуємо простір Соболєва $\mathcal{B}_{2,k}(\mathbb{R}^{n})=H^{(s)}(\mathbb{R}^{n})$.

Простори $\mathcal{B}_{2,k}(\mathbb{R}^{n})$ посідають центральне місце серед просторів узагальненої гладкості \cite{FarkasLeopold06, HaroskeMoura08, KalyabinLizorkin87}. Вони є об'єктом різних досліджень, частину з яких  виконано в останні десять років (див., наприклад, \cite{BesoyCobos18, DominguezTikhonov23, HaroskeLeopoldMouraSkrzypczak23, LoosveldtNicolay19, MouraNevesSchneider14, NevesOpic20}). Чимало з них пов'язано з теорією еліптичних операторів і еліптичних крайових задач та викладено в монографії \cite{MikhailetsMurach14} і оглядах \cite{MikhailetsMurach09OperatorTheory191, MikhailetsMurach12BJMA2, MikhailetsMurach11NaukVisnChernivtci} перших двох авторів. Даний огляд містить результати, які отримано в останні десять років і доповнюють дослідження, викладені у цих працях. Їх мотивація більш детально  висвітлена у вступах до \cite{MikhailetsMurach12BJMA2, MikhailetsMurach14}. Доцільність та перспективи таких досліджень обговорені у монографії Ґ.~Трібеля \cite[сс. 57--60]{Triebel10}.

Зауважимо, що в останній час простори узагальненої гладкості знайшли важливі застосування до різних параболічних диференціальних рівнянь \cite{DyachenkoLos22JEPE1, DyachenkoLos23UMJ8, Los17UMJ3, LosMikhailetsMurach17CPAA1, LosMikhailetsMurach21CPAA10, LosMurach17OpenMath, MikuleviciusPhonsom19, MikuleviciusPhonsom21, LosMikhailetsMurach21monograph}. Зважаючи на це, авторами готується огляд застосувань просторів Хермандера до параболічних початково-крайових задач.

Стаття складається з 12 розділів. Перший з них~--- вступ. Розділ~\ref{survey-sect2} присвячений RO-змінним функціям на нескінченності, які служать показником узагальненої гладкості для гільбертових просторів Хермандера, розглянутих у розд.~\ref{survey-sect3}. Ці простори утворюють розширену соболєвську шкалу на $\mathbb{R}^{n}$ та евклідових областях. Її інтерполяційні властивості обговорено у розд.~\ref{survey-sect4}. У розд.~\ref{survey-sect5} розглянуто таку шкалу на компактному $C^{\infty}$-многовиді без краю.

Розділи \ref{survey-sect6}--\ref{survey-sect8} присвячені застосуванням розширеної соболєвської шкали до лінійних еліптичних диференціальних систем. У розд.~\ref{survey-sect6} розглянуто системи в $\mathbb{R}^{n}$, рівномірно еліптичні за Дуглісом--Ніренбергом. Він містить глобальну апріорну оцінку їх розв'язків, достатні умови їх узагальненої та класичної гладкості та умови фредгольмовості відповідного матричного диференціального оператора на підходящих парах просторів Хермандера. У розд.~\ref{survey-sect7} обговорені версії цих результатів для еліптичних за Дуглісом--Ніренбергом диференціальних систем, заданих на компактному $C^{\infty}$-многовиді без краю. Еліптичні системи з параметром на такому многовиді розглянуті в розд.~\ref{survey-sect8}. Оператори, породжені цими системами встановлюють ізоморфізми на парах просторів Хермандера для великих за модулем значень комплексного параметра з кута еліптичності системи, а розв'язки допускають двобічну апріорну оцінку зі сталими, не залежними від цих значень.

Розділи \ref{survey-sect9}--\ref{survey-sect12} присвячені застосуванням розширеної соболєвської шкали до лінійних еліптичних диференціальних крайових задач. У розд.~\ref{survey-sect9} розглянуто регулярні еліптичні крайові задачі, задані в обмеженій евклідовій області класу $C^{\infty}$. Зазначено, що ці задачі є фредгольмовими на підходящих парах просторів Хермандера і наведено достатні умови узагальненої та класичної гладкості розв'язків, а також їх локальну (аж до частини межі області) апріорну оцінку у просторах Хермандера. У~розд.~\ref{survey-sect10} розглянуто еліптичні крайові задачі з параметром у цих просторах. Він містить версії результатів розд.~\ref{survey-sect8} для цих задач. Властивості регулярних еліптичних задач з крайовими даними як завгодно низької регулярності обговорено у розд.~\ref{survey-sect11}. Вони доповнюють результати розд.~\ref{survey-sect9}. Окремо розглянуто випадок однорідного еліптичного диференціального рівняння. Останній розділ~\ref{survey-sect12} присвячений еліптичним крайовим задачам з додатковими невідомими функціями у крайових умовах, порядки яких можуть бути як завгодно великими. Він містить теореми про фредгольмовість цих задач на парах просторів Хермандера, умови локальної регулярності їх розв'язків та відповідні їх апріорні оцінки.

Сформульовані теореми порівняно з близькими результатами, отриманими для вужчих або інших класів просторів узагальненої гладкості, еліптичних рівнянь і еліптичних крайових задач.

\section{RO-змінні функції}\label{survey-sect2}

Ці функції служать показником гладкості для гільбертових просторів Хермандера, які утворюють розширену соболєвську шкалу.

За означенням, клас $\mathrm{RO}$ складається з усіх вимірних за Борелем функцій $\alpha:\nobreak[1,\infty)\rightarrow(0,\infty)$, для яких існують числа $b>1$ і $c\geq1$ такі, що $c^{-1}\leq\alpha(\lambda t)/\alpha(t)\leq c$ для усіх $t\geq1$ і $\lambda\in[1,b]$ (числа $b$ і $c$ можуть залежати від $\alpha$). Такі функції називають RO- (або OR-) \textit{змінними на нескінченності}.

Клас $\mathrm{RO}$ введений В.~Г.~Авакумовичем \cite{Avakumovic36}, досить повно вивчений і має різні застосування, зокрема у теоремах тауберового типу \cite{BinghamGoldieTeugels89, BuldyginIndlekoferKlesovSteinebach18, Seneta76}.

Кожна функція $\alpha\in\mathrm{RO}$ обмежена і відокремлена від нуля на будь-якому відрізку $[1,q]$, де $1<q<\infty$ \cite[лема~A.1]{Seneta76}.
Клас $\mathrm{RO}$ має такий інтегральний опис \cite[теорема~A.1]{Seneta76}.

\begin{theorem}
Функція $\alpha:\nobreak[1,\infty)\rightarrow(0,\infty)$ належить до $\mathrm{RO}$ тоді і тільки тоді, коли
\begin{equation}\label{description-RO}
\alpha(t)=\exp\Biggl(\beta(t)+
\int\limits_{1}^{\:t}\frac{\gamma(\tau)}{\tau}\;d\tau\Biggr)\quad\mbox{при}
\quad t\geq1
\end{equation}
для деяких вимірних за Борелем обмежених функцій $\beta,\gamma:\nobreak[1,\infty)\rightarrow\mathbb{R}$.
\end{theorem}

Важливу роль для функцій класу $\mathrm{RO}$  відіграють їх індекси Матушевської~\cite{Matuszewska64}. Для означення цих індексів знадобиться така властивість класу RO \cite[теорема~A.2]{Seneta76}.

\begin{theorem}\label{prop-bounds}
Для будь-якої функції $\alpha\in\mathrm{RO}$ існують дійсні числа
$s_{0}$ і $s_{1}$, причому $s_{0}\leq s_{1}$, і додатні числа $c_{0}$ і $c_{1}$ такі, що
\begin{equation}\label{21.1}
c_{0}\lambda^{s_{0}}\leq\frac{\alpha(\lambda t)}{\alpha(t)}\leq c_{1}\lambda^{s_{1}}\;\;\;\mbox{для усіх}\;\;\;t\geq1\;\;\; \mbox{і}\;\;\;\lambda\geq1.
\end{equation}
\end{theorem}

Нижній і верхній індекси Матушевської функції $\alpha\in\mathrm{RO}$ означаються відповідно за формулами
\begin{gather}\label{deq211}
\sigma_{0}(\alpha):= \sup\{s_{0}\in\mathbb{R}:\,\mbox{виконується
ліва нерівність в \eqref{21.1}}\},\\
\sigma_{1}(\alpha):=\inf\{s_{1}\in\mathbb{R}:\,\mbox{виконується
права нерівність в \eqref{21.1}}\}\label{deq212}
\end{gather}
\cite[п.~2.1.2]{BinghamGoldieTeugels89}; тут $-\infty<\sigma_{0}(\alpha)\leq\sigma_{1}(\alpha)<\infty$.

З теореми \ref{prop-bounds} випливає, що кожна функція $\alpha\in\mathrm{RO}$ є міжстепеневою, тобто
\begin{equation}\label{21.2}
(s_{0}<\sigma_{0}(\alpha),\;\sigma_{1}(\alpha)<s_{1})\;\Longrightarrow\;
(c_{0}'\lambda^{s_{0}}\leq\alpha(\lambda)\leq c_{1}'\lambda^{s_{1}}\;\mbox{для усіх}\;\lambda\geq1),
\end{equation}
де додатні числа $c_{0}'$ і $c_{1}'$ не залежать від $\lambda$.

Вкажемо відомі \cite[вправи до розд.~1]{Seneta76} приклади функцій з класу $\mathrm{RO}$.

%%\medskip

\begin{example}\label{example2.1}
До класу $\mathrm{RO}$ належить кожна неперервна функція $\alpha:[1,\infty)\rightarrow(0,\infty)$ така, що
\begin{equation*}%%\label{21.3}
\alpha(t)=t^{s}(\ln t)^{r_{1}}(\ln\ln
t)^{r_{2}}\ldots(\underbrace{\ln\ldots\ln}_{k\;\mbox{\small
разів}} t)^{r_{k}}\quad\mbox{для}\quad t\gg1,
\end{equation*}
де числа $k\in\mathbb{N}$ і $s,r_{1},\ldots,r_{k}\in\mathbb{R}$ вибрані довільно. Для неї $\sigma_{0}(\alpha)=\sigma_{1}(\alpha)=s$.
\end{example}

\begin{example}\label{1ex4}
До класу $\mathrm{RO}$ належить функція $\alpha(t):=t^{s}e^{(\ln t)^{r}}$  аргументу $t\geq1$, де числа $s\in\mathbb{R}$ і $r\in(0,1)$ вибрані довільно. Для неї також $\sigma_{0}(\alpha)=\sigma_{1}(\alpha)=s$.
\end{example}

Ці функції є еталонним прикладом функцій, правильно змінних за Й.~Караматою \cite{Karamata30a} на нескінченності.

Вимірна за Борелем функція $\alpha:\nobreak[1,\infty)\rightarrow(0,\infty)$ називається \textit{правильно змінною} за Караматою на нескінченності, якщо існує число $s\in\mathbb{R}$ таке, що $\alpha(\lambda t)/\alpha(t)\to\lambda^{s}$ при $t\rightarrow\infty$ для кожного $\lambda>0$. Число $s$ називається \textit{порядком змінення} функції $\alpha$ на нескінченності. У випадку  $s=0$ ця функція називається \textit{повільно змінною} на нескінченності за Караматою.

Існують функції $\alpha\in\mathrm{RO}$ такі, що
$\sigma_0(\alpha)=\sigma_1(\alpha)$, але $\alpha$ не
еквівалентна в околі нескінченності жодній функції $\eta>\nobreak0$, правильно змінній на нескінченності. (Як звичайно, додатні функції $\alpha$ і $\eta$ називаємо еквівалентними в околі нескінченності, якщо обидві функції $\alpha(t)/\eta(t)$ і $\eta(t)/\alpha(t)$ обмежені на деякому промені $[r,\infty)$, де $r\gg1$.)

\begin{example}
Прикладом такої функції $\alpha\in\mathrm{RO}$ служить функція $\alpha(t):=e^{h(\ln t)}$ аргументу $t\geq1$, де $h$ означено за формулами: $h(x):=0$ при $x\in[0,1]$ та $h(x):=h(2^j)+(x-2^j)^{1/2}$ при $x\in[2^j,2^{j+1}]$ для кожного цілого $j\geq0$; для неї $\sigma_0(\alpha)=\sigma_1(\alpha)=0$ \cite[твердження~2.2.8]{BinghamGoldieTeugels89}.
\end{example}

Наведемо приклади функції класу $\mathrm{RO}$ з різними індексами Матушевської.

\begin{example}
Виберемо довільно числа $\theta\in\mathbb{R}$, $\delta>0$, $r\in(0,1]$ і
означимо функцію $\alpha$ за формулами $\alpha(t):=t^{\theta+\delta\sin(\ln\ln t)^{r}}$, якщо $t>e$, і $\alpha(t):=t^{\theta}$, якщо $1\leq t\leq e$. Тоді $\alpha\in\mathrm{RO}$, причому $\sigma_{0}(\alpha)=\theta-\delta$ і $\sigma_{1}(\alpha)=\theta+\delta$ у випадку $0<r<1$, але $\sigma_{0}(\alpha)=\theta-\sqrt{2}\delta$ і $\sigma_{1}(\alpha)=\theta+\sqrt{2}\delta$ у випадку $r=1$. Якщо $r>1$, то $\alpha\notin\mathrm{RO}$ \cite[п.~3]{AnopDenkMurach21}.
\end{example}

\begin{example}
Виберемо довільно дійсні числа $r<s$ і зростаючу послідовність $(\theta_{k})_{k=1}^{\infty}$ таку, що $\theta_{1}=1$ та $\theta_{k}\to\infty$ при $k\to\infty$. Покладемо $\pi_{k}:=\theta_{1}\ldots\theta_{k}$ для кожного $k\in\mathbb{N}$ і означимо функцію $\alpha\in\mathrm{RO}$ за формулою \eqref{description-RO} де $\beta(t)\equiv0$, а $\gamma(\tau):=r$, якщо $\tau\in[\pi_{2j-1},\pi_{2j}]$ для деякого $j\in\mathbb{N}$, та $\gamma(\tau):=s$ у противному разі. Для цієї функції  $\sigma_{0}(\alpha)=r$ і $\sigma_{1}(\alpha)=s$ \cite[п.~3]{AnopChepurukhinaMurach21Axioms}.
\end{example}

\section{Розширена соболєвська шкала}\label{survey-sect3}

Ця шкала складається з функціональних просторів Хермандера $H^{\alpha}$, для яких показник гладкості $\alpha$ пробігає клас $\mathrm{RO}$. Усі лінійні простори вважаємо комплексними, а функції та розподіли~--- комплекснозначними (якщо інше не зазначено окремо). Розглянемо спочатку простори на $\mathbb{R}^{n}$, де $\nobreak{n\in\mathbb{N}}$. Нехай $\alpha\in\mathrm{RO}$.

За означенням, лінійний простір $H^{\alpha}(\mathbb{R}^{n})$ складається з усіх розподілів $w\in\mathcal{S}'(\mathbb{R}^{n})$ таких, що їх перетворення Фур'є $\widehat{w}$ локально інтегровне за Лебегом на
$\mathbb{R}^{n}$ і задовольняє умову
\begin{equation*}
\int_{\mathbb{R}^{n}}
\alpha^2(\langle\xi\rangle)\,|\widehat{w}(\xi)|^2\,d\xi<\infty.
\end{equation*}
У просторі  $H^{\alpha}(\mathbb{R}^{n})$ означені скалярний
добуток і норма за формулами
\begin{equation*}%%\label{21.4}
(w_1,w_2)_{\alpha,\mathbb{R}^{n}}:=
\int_{\mathbb{R}^{n}}\alpha^2(\langle\xi\rangle)\,
\widehat{w_1}(\xi)\,\overline{\widehat{w_2}(\xi)}\,d\xi
\quad\mbox{і}\quad
\|w\|_{\alpha,\mathbb{R}^{n}}:=
(w,w)_{\alpha,\mathbb{R}^{n}}^{1/2}.
\end{equation*}

Тут, як звичайно, $\langle\xi\rangle:=(1+|\xi|^{2})^{1/2}$~--
згладжений модуль вектора $\xi\in\mathbb{R}^{n}$, а $\mathcal{S}'(\mathbb{R}^{n})$~-- лінійний топологічний простір Л.~Шварца повільно зростаючих розподілів в $\mathbb{R}^{n}$. Всі розподіли трактуються як \emph{анти}лінійні функціонали на відповідних просторах основних функцій. Отже, $\mathcal{S}'(\mathbb{R}^{n})$ є антидуальним простором до лінійного топологічного простору $\mathcal{S}(\mathbb{R}^{n})$ нескінченно гладких швидко спадаючих функцій на $\mathbb{R}^{n}$, причому дуальність розглядається відносно скалярного добутку у гільбертовому просторі $L_{2}(\mathbb{R}^{n})$ функцій, квадратично інтегровних на $\mathbb{R}^{n}$.

Простір $H^{\alpha}(\mathbb{R}^{n})$ повний (тобто гільбертів) і сепарабельний. Виконуються неперервні і щільні вкладення $\mathcal{S}(\mathbb{R}^{n})\hookrightarrow
H^{\alpha}(\mathbb{R}^{n})\hookrightarrow\mathcal{S}'(\mathbb{R}^{n})$. У просторі $H^{\alpha}(\mathbb{R}^{n})$ також щільна множина $C_{0}^{\infty}(\mathbb{R}^{n})$ усіх нескінченно диференційовних на $\mathbb{R}^{n}$ функцій з компактним носієм

Простір $H^{\alpha}(\mathbb{R}^{n})$ є гільбертів ізотропний
випадок просторів $\mathcal{B}_{p,k}(\mathbb{R}^{n})$, введених і досліджених Л.~Хермандером \cite[п.~2.2]{Hermander63} (див також його монографію \cite[п.~10.1]{Hermander05v2}). А саме,
$H^{\alpha}(\mathbb{R}^{n})=\mathcal{B}_{p,k}(\mathbb{R}^{n})$, якщо $p=2$ і $k(\xi)=\alpha(\langle\xi\rangle)$ для усіх $\xi\in\mathbb{R}^{n}$.
У гільбертовому випадку $p=2$ простори Хермандера
збігаються з просторами, введеними і дослідженими Л.~Р.~Волевичем і
Б.~П.~Панеяхом \cite[\S~2]{VolevichPaneah65}.

У випадку степеневої функції $\alpha(t)\equiv t^{s}$, простір
$H^{\alpha}(\mathbb{R}^{n})$ стає гільбертовим простором
Соболєва $H^{(s)}(\mathbb{R}^{n})$ порядку $s\in\mathbb{R}$. У загальному випадку виконуються щільні неперервні вкладення
\begin{equation}\label{21.6}
H^{(s_1)}(\mathbb{R}^{n})\hookrightarrow
H^{\alpha}(\mathbb{R}^{n})\hookrightarrow
H^{(s_0)}(\mathbb{R}^{n}),
\quad\mbox{якщо}\quad
s_{0}<\sigma_{0}(\alpha)\quad\mbox{і}\quad\sigma_{1}(\alpha)<s_{1}.
\end{equation}
Це випливає з властивості \eqref{21.2}.

Клас гільбертових функціональних просторів $\{H^{\alpha}(\mathbb{R}^{n}):\alpha\in\mathrm{RO}\}$
називаємо \textit{розширеною соболєвською шкалою} на $\mathbb{R}^{n}$.
Він  виділений і досліджений в \cite{MikhailetsMurach13UMJ3, MikhailetsMurach15ResMath1, MikhailetsMurach08Collection1} і
\cite[п.~2.4.2]{MikhailetsMurach14}, а його назва запропонована в \cite[п.~2]{MikhailetsMurach13UMJ3}.

З властивостей просторів Хермандера \cite[п.~2.2]{Hermander63} випливають такі властивості розширеної соболєвської шкали (\cite[Твердження~2]{MikhailetsMurach13UMJ3} або \cite[Твердження~2.6]{MikhailetsMurach14}).

\begin{theorem}\label{th-properties-H}
Нехай $\alpha,\alpha_{1}\in\mathrm{RO}$. Тоді:
\begin{itemize}
  \item[(i)] Функція $\alpha(t)/\alpha_{1}(t)$ обмежена в околі нескінченності тоді і тільки тоді, коли $H^{\alpha_{1}}(\mathbb{R}^{n})\hookrightarrow H^{\alpha}(\mathbb{R}^{n})$. Це вкладення неперервне і щільне.
  \item[(ii)] Простори $H^{\alpha}(\mathbb{R}^{n})$ і
$H^{1/\alpha}(\mathbb{R}^{n})$ взаємно дуальні відносно розширення за неперервністю скалярного добутку в $L_{2}(\mathbb{R}^{n})$.
  \item[(iii)] Для кожного цілого числа $\lambda\geq0$ нерівність
\begin{equation*}
\int_{1}^{\,\infty}\;t^{2\lambda+n-1}\alpha^{-2}(t)\,dt<\infty
\end{equation*}
рівносильна вкладенню $H^{\alpha}(\mathbb{R}^{n})\hookrightarrow
C^{\lambda}_{\mathrm{b}}(\mathbb{R}^{n})$. Це вкладення неперервне.
\end{itemize}
\end{theorem}

Стосовно твердження~(іі) зауважимо, що
$\alpha\in\mathrm{RO}\Leftrightarrow1/\alpha\in\mathrm{RO}$. Як звичайно, у твердженні~(ііі) через $C^{\lambda}_{\mathrm{b}}(\mathbb{R}^{n})$ позначено банахів простір усіх функцій, які задані на $\mathbb{R}^{n}$ і мають неперервні і обмежені на $\mathbb{R}^{n}$ частинні похідні до порядку $\lambda$ включно.

Нехай $\Omega$~--- довільна непорожня відкрита множина в $\mathbb{R}^{n}$. Розширена соболєвська шкала на $\Omega$~--- це клас  $\{H^{\alpha}(\Omega):\alpha\in\mathrm{RO}\}$ просторів, які означаються у такий спосіб: лінійний простір $H^{\alpha}(\Omega)$ складається зі звужень в $\Omega$ усіх розподілів $w\in H^{\alpha}(\mathbb{R}^{n})$ і наділений нормою
\begin{equation*}%%\label{21.9}
\|u\|_{\alpha,\Omega}:=
\inf\bigl\{\|w\|_{\alpha,\mathbb{R}^{n}}:\,
w\in H^{\alpha}(\mathbb{R}^{n}),\;w=u\;\,\mbox{на}\;\,\Omega\bigr\},
\end{equation*}
де $u\in H^{\alpha}(\Omega)$. Простір $H^{\alpha}(\Omega)$ гільбертів і сепарабельний відносно цієї норми, бо він~--- факторпростір гільбертового сепарабельного простору $H^{\alpha}(\mathbb{R}^{n})$ за підпростором
\begin{equation}\label{21.10}
\bigl\{\omega\in H^{\alpha}(\mathbb{R}^{n}):
\mathrm{supp}\,\omega\subset\mathbb{R}^{n}\setminus\Omega\bigr\}.
\end{equation}
Норма у просторі $H^{\alpha}(\Omega)$ породжена скалярним добутком
\begin{equation*}%%\label{21.11}
(u_1,u_2)_{\alpha,\Omega}:=
(w_1-\Pi w_1,w_2-\Pi w_2)_{\alpha,\mathbb{R}^{n}},
\end{equation*}
де $u_j\in H^\alpha(\Omega)$, $w_j\in
H^{\alpha}(\mathbb{R}^{n})$, $u_j=w_j$ на $\Omega$ для кожного $j\in\{1,2\}$, а $\Pi$~--- ортопроектор простору $H^{\alpha}(\mathbb{R}^{n})$ на підпростір \eqref{21.10}.

Відмітимо, що $H^{\alpha}(\Omega)$~--- окремий ізотропний випадок гільбертових функціональних просторів, які ввели і дослідили
Л.~Р.~Волевич і Б.~П.~Панеях \cite[\S~3]{VolevichPaneah65}. Якщо $\alpha(t)\equiv t^{s}$ для деякого $s\in\mathbb{R}$, то $H^{\alpha}(\Omega)=H^{(s)}(\Omega)$~--- гільбертів простір Соболєва порядку $s$, заданий на множині $\Omega$. При цьому приймаємо досить поширене означення просторів Соболєва на евклідовій області, однакове як для додатного, так і для від'ємного порядку (див., наприклад, монографію Г.~Трібеля \cite[п.~4.2.1]{Triebel95}).

З наведеного означення випливає, що твердження (i), (iii)  теореми~\ref{th-properties-H} і  властивість \eqref{21.6} зберігають силу, якщо у них перейти до звужень функціональних просторів на $\Omega$. Якщо множина $\Omega$ обмежена, то для неї аналог твердження~(i) цієї теореми має таке уточнення: вкладення $H^{\alpha_{1}}(\Omega)\hookrightarrow H^{\alpha}(\Omega)$ компактне тоді і тільки тоді, коли $\alpha(t)/\alpha_{1}(t)\to\nobreak0$ при $t\to\infty$; це випливає з \cite[теорема~2.2.3]{Hermander63}. Тому є компактними вкладення в аналогах твердження~(iii) і властивості \eqref{21.6} для обмеженої відкритої множини~$\Omega$.

Якщо $\sigma_{0}(\alpha)=\sigma_{1}(\alpha)=:s$, то $\alpha(t)\equiv t^{s}\varphi(t)$, де $\varphi\in\mathrm{RO}$ і $\sigma_{0}(\varphi)=\sigma_{1}(\varphi)=0$. У цьому випадку простір $H^{\alpha}(\Omega)$ зручно позначити через $H^{s,\varphi}(\Omega)$, вказуючи основний числовий показник гладкості $s$ і додатковий функціональний показник гладкості $\varphi$. Якщо функція $\varphi$ є повільно змінною на нескінченності за Караматою, то простори $H^{s,\varphi}(\Omega)$ утворюють \textit{уточнену соболєвську шкалу} на $\Omega$, введену і досліджену в \cite{MikhailetsMurach05UMJ5, MikhailetsMurach06UMJ3, MikhailetsMurach06Dop6}. У~статтях В.~С.~Ільківа, Н.~І.~Страп, І.~І.~Волянської \cite{IlkivStrap15, IlkivStrapVolyanska20} уведено аналог цієї шкали для просторів многочленів від $2n$ комплексних змінних $z_{1},\ldots,z_{n}$ і $z_{1}^{-1},\ldots,z_{n}^{-1}$ та наведено її застосування до нелокальних крайових задач для диференціально-операторних рівнянь.

\section{Інтерполяційні властивості шкали}\label{survey-sect4}

Вони пов'язані з методом квадратичної інтерполяції (з функціональним параметром) пар гільбертових просторів. Він уперше з'явився у статті Ч.~Фояша і Ж.-Л.~Ліонса \cite[с.~278]{FoiasLions61} і є природнім узагальненням класичного інтерполяційного методу Ж.-Л.~Ліонса \cite{Lions58} і C.~Г.~Крейна \cite{Krein60a} на випадок, коли замість числового параметра інтерполяції береться досить загальна функція. Цей метод інтерполяції та пов'язані з ним гільбертові шкали досліджено у працях \cite{Ameur04, Ameur19, Donoghue67, Fan11, KreinPetunin66, MikhailetsMurach06UMJ2, MikhailetsMurach08MFAT1, Ovchinnikov84, Pustylnik82} та викладено у монографіях \cite[розд.~1]{LionsMagenes72},
\cite[п.~1.1]{MikhailetsMurach14} і \cite[розд. 15,~30]{Simon19}.
Наведемо означення цього методу, слідуючи в основному  \cite[п.~1.1]{MikhailetsMurach14} і обмежуючись сепарабельними гільбертовими просторами.

Нехай упорядкована пара $X:=[X_{0},X_{1}]$ гільбертових просторів є \textit{регулярною}, тобто $X_{1}$ є щільним лінійним многовидом у $X_{0}$ і вкладення $X_{1}\subset X_{0}$ неперервне. Тоді існує додатно визначений самоспряжений оператор $J$ у гільбертовому просторі $X_{0}$ з областю визначення $X_{1}$ такий, що $\|Jw\|_{X_{0}}=\|w\|_{X_{1}}$ для довільного $w\in X_{1}$ \cite[розд.~1, п.~2.1]{LionsMagenes72}. Оператор $J$ визначається за парою $X$ однозначно і називається \emph{породжуючим} оператором для неї.

Нехай вимірна за Борелем функція $\psi:(0,\infty)\rightarrow(0,\infty)$  обмежена на кожному відрізку $[a,b]$ і відокремлена від нуля на кожній множині $[r,\infty)$, де $0<a<b<\infty$ і $r>0$. За допомогою спектрального розкладу самоспряженого оператора $J$ означено самоспряжений оператор $\psi(J)$ у просторі $X_{0}$. Позначимо через $[X_{0},X_{1}]_{\psi}$ або $X_{\psi}$ область визначення оператора $\psi(J)$, наділену скалярним добутком і нормою
\begin{equation*}
(w_1,w_2)_{X_{\psi}}:=(\psi(J)w_1,\psi(J)w_2)_{X_{0}}
\quad\mbox{і}\quad \|w\|_{X_{\psi}}:=\|\psi(J)w\|_{X_{0}}.
\end{equation*}
Простір $X_{\psi}$ гільбертів і сепарабельний.

Функція $\psi$ називається \emph{інтерполяційним
параметром}, якщо для будь-яких регулярних пар $X=[X_{0},X_{1}]$ і
$Y=[Y_{0},Y_{1}]$ гільбертових просторів та для довільного лінійного
відображення $T$, заданого на $X_{0}$, виконується така властивість:
якщо для кожного індексу $j\in\{0,1\}$ звуження відображення $T$ на
простір $X_{j}$ є обмеженим оператором $T:X_{j}\rightarrow Y_{j}$,
то і звуження відображення $T$ на простір $X_{\psi}$ є обмеженим
оператором $T:X_{\psi}\rightarrow Y_{\psi}$. Якщо функція $\psi$ є інтерполяційним параметром, то говоримо, що простір $X_{\psi}$ отримується методом квадратичної інтерполяції з функціональним параметром $\psi$ пари $X$, а оператор $T:X_{\psi}\rightarrow Y_{\psi}$ є результатом інтерполяції операторів $T:X_{j}\rightarrow Y_{j}$, де $j\in\{0,1\}$. Тоді виконуються неперервні та щільні вкладення $X_{1}\hookrightarrow X_{\psi}\hookrightarrow X_{0}$ \cite[теорема~1.1]{MikhailetsMurach14}.

З  теореми Ж.~Петре \cite{Peetre66, Peetre68} про опис усіх
інтерполяційних функцій додатного порядку (див. також монографію
\cite[теорема~5.4.4]{BerghLefstrem76}) випливає такий результат: функція $\psi$ є інтерполяційним параметром тоді і тільки тоді, коли вона псевдоугнута в околі нескінченності, тобто еквівалентна там деякій додатній угнутій функції \cite[теорема~1.9]{MikhailetsMurach14}.

Припустимо, що $\Omega$ є або $\mathbb{R}^{n}$, або відкритим підпростором в $\mathbb{R}^{n}$, або обмеженою ліпшіцевою областю в $\mathbb{R}^{n}$. Розширена соболєвська шкала на $\Omega$ має такі важливі інтерполяційні властивості: отримується квадратичною інтерполяцією пар соболєвських просторів, замкнута відносно квадратичної інтерполяції та складається з усіх інтерполяційних гільбертових просторів щодо пар просторів Соболєва. Сформулюємо відповідні результати \cite[теореми 5.1, 5.2, 2.4]{MikhailetsMurach15ResMath1}.

\begin{theorem}\label{th-int-Sobolev}
Нехай задано функцію $\alpha\in\mathrm{RO}$ і числа $s_{0},s_{1}\in\mathbb{R}$ такі, що $s_{0}<\sigma_{0}(\alpha)$ і $\sigma_{1}(\alpha)<s_{1}$. Покладемо
\begin{equation*}%%\label{21.25}
\psi(t):=
\begin{cases}
\;t^{{-s_0}/{(s_1-s_0)}}\,
\alpha(t^{1/{(s_1-s_0)}}),&\text{якщо}\quad t\geq1, \\
\;\alpha(1),&\text{якщо}\quad 0<t<1.
\end{cases}
\end{equation*}
Тоді функція $\psi$~--- інтерполяційний параметр і
\begin{equation*}%%\label{21.26}
H^{\alpha}(\Omega)=
\bigl[H^{(s_{0})}(\Omega),H^{(s_{1})}(\Omega)\bigr]_{\psi}
\end{equation*}
з еквівалентністю норм (з рівністю норм, якщо $\Omega=\mathbb{R}^{n}$).
\end{theorem}

\begin{theorem}\label{th-int-Hormander}
Нехай $\alpha_{0},\alpha_{1}\in\mathrm{RO}$, причому функція $\alpha_{0}/\alpha_{1}$ обмежена в околі нескінченності. Припустимо, що функція $\psi$~--- інтерполяційний параметр, і покладемо
\begin{equation*}
\alpha(t):=\alpha_{0}(t)\,\psi
\left(\frac{\alpha_{1}(t)}{\alpha_{0}(t)}\right)\quad\mbox{для довільного}\quad t\geq1.
\end{equation*}
Тоді $\alpha\in\mathrm{RO}$ і
\begin{equation*}
[H^{\alpha_{0}}(\Omega),H^{\alpha_{1}}(\Omega)]_{\psi}= H^{\alpha}(\Omega)
\end{equation*}
з еквівалентністю норм (з рівністю норм, якщо $\Omega=\mathbb{R}^{n}$).
\end{theorem}

Стосовно наступної властивості нагадаємо, що гільбертів простір $H$ називається \textit{інтерполяційним} для регулярної пари $X=[X_0,X_1]$ гільбертових просторів, якщо: а)~виконуються неперервні вкладення $X_1\hookrightarrow H\hookrightarrow X_0$; б)~з того, що лінійне відображення $T$ є обмеженим оператором на $X_0$ і на $X_1$, випливає його обмеженість на $H$. З результату В.~І.~Овчинникова \cite[теорема~14.4.1]{Ovchinnikov84} випливає, що гільбертів простір $H$ є інтерполяційним для вказаної пари $X=[X_0,X_1]$ тоді і лише тоді, коли $H=X_{\psi}$ з еквівалентністю норм для деякого інтерполяційного параметра~$\psi$.

\begin{theorem}\label{th-int-space}
Нехай задано дійсні числа $s_{0}<s_{1}$. Гільбертів простір $H$ є інтерполяційним для пари просторів Соболєва $[H^{(s_{0})}(\Omega),H^{(s_{1})}(\Omega)]$
тоді і лише тоді, коли $H=H^{\alpha}(\Omega)$ з еквівалентністю норм для деякого функціонального параметра $\alpha\in\mathrm{RO}$, який задовольняє умову~\eqref{21.1}.
\end{theorem}

Зауважимо, що \eqref{21.1} еквівалентна такий парі умов, поданих у термінах індексів Матушевської функції $\alpha\in\mathrm{RO}$:
\begin{enumerate}
\item [$\mathrm{(i)}$] $s_{0}\leq\sigma_{0}(\alpha)$ та, крім того,
$s_{0}<\sigma_{0}(\alpha)$ якщо не досягається супремум в \eqref{deq211};
\item [$\mathrm{(ii)}$] $\sigma_{1}(\alpha)\leq s_{1}$ та, крім того, $\sigma_{1}(\alpha)<s_{1}$ якщо не досягається інфімум в $\eqref{deq212}$.
\end{enumerate}

З теорем \ref{th-int-Sobolev} і \ref{th-int-Hormander} випливають інтерполяційні нерівності для норм у просторах, які утворюють розширену соболєвську шкалу. Це дозволяє інтерпретувати її як змінну гільбертову шкалу, досліджену, наприклад, в \cite{Hegland10, MatheTautenhahn06}.

У новітній статті О.~Мілатовича \cite{Milatovic24} уведена і досліджена розширена соболєвська шкала на решітці $\mathbb{Z}^{n}$. Для цієї шкали доведено версії інтерполяційних теорем \ref{th-int-Sobolev}--\ref{th-int-space} і теореми \ref{th-properties-H}.

При інтерполяції просторів успадковується не лише обмеженість лінійних операторів, а і їх фредгольмовість за певних умов. Це важливо у теорії еліптичних диференціальних операторів, які є фредгольмовими на відповідних парах соболєвських просторів.

Обмежений лінійний оператор
$T:E_{1}\rightarrow E_{2}$, де $E_{1}$ і $E_{2}$~--- банахові простори, називається \textit{фредгольмовим} (або \textit{нетеровим}), якщо
його ядро $\ker T:=\{w\in E_{1}:Tw=0\}$ і коядро $\mathrm{coker}\,T:=E_{2}/T(E_{1})$ скінченновимірні. Якщо оператор $T$ фредгольмів, то його область значень $T(E_{1})$ замкнена у просторі $E_{2}$ \cite[лема~19.1.1]{Hermander07v2}, а \textit{індекс} $\mathrm{ind}\,T:=\dim\ker T-\dim\mathrm{coker}\,T$ скінченний.

Нехай $X=[X_0,X_1]$ і $Y=[Y_0,Y_1]$ є регулярні пари гільбертових
просторів. Нехай, окрім того, на $X_0$ задане лінійне відображення
$T$ таке, що його звуження на простори $X_j$, де $j=0,1$, є фредгольмовими обмеженими операторами $T:X_j\rightarrow Y_j$, які
мають спільне ядро і однаковий індекс. Тоді для довільного
інтерполяційного параметра $\psi$ обмежений оператор
$T:X_\psi\rightarrow Y_\psi$ фредгольмів з тим же ядром і індексом, а
його область значень дорівнює $Y_\psi\cap T(X_0)$.

Цей результат \cite[теорема~1.7]{MikhailetsMurach14} правильний для будь-якого методу інтерполяції пар банахових просторів \cite[теорема~5.1]{Geymonat65}.

\section{Шкала на компактному многовиді}\label{survey-sect5}

Нехай $\Gamma$~--- нескінченно гладкий компактний многовид без краю (тобто замкнений многовид) вимірності $n\in\mathbb{N}$ і наділений деякою щільністю $dx$ класу $C^{\infty}$. Виберемо довільно скінченний атлас із $C^{\infty}$-структури на многовиді $\Gamma$. Нехай він утворений локальними картами $\pi_j: \mathbb{R}^{n}\leftrightarrow \Gamma_{j}$, де $j=1,\ldots,p$. Тут відкриті множини $\Gamma_{1},\ldots,\Gamma_{p}$ складають покриття многовиду $\Gamma$. Крім того, виберемо довільно функції $\chi_j\in C^{\infty}(\Gamma)$, де $j=1,\ldots,p$, які утворюють розбиття одиниці на $\Gamma$, підпорядковане умові $\mathrm{supp}\,\chi_j\subset \Gamma_j$.

Нехай $\alpha\in\mathrm{RO}$. За означенням, лінійний простір $H^{\alpha}(\Gamma)$ складається з усіх розподілів $h\in\mathcal{D}'(\Gamma)$, які в локальних картах дають елементи простору $H^{\alpha}(\mathbb{R}^{n})$, тобто  $(\chi_{j}h)\circ\pi_{j}\in H^{\alpha}(\mathbb{R}^{n})$ для кожного номера $j\in\{1,\ldots,p\}$. Скалярний добуток і норма у просторі $H^{\alpha}(\Gamma)$ означені за формулами
\begin{equation*}
(h_{1},h_{2})_{\alpha,\Gamma}:=
\sum_{j=1}^{p}\,((\chi_{j}h_{1})\circ\pi_{j},
(\chi_{j}\,h_{2})\circ\pi_{j})_{\alpha,\mathbb{R}^{n}}
\quad\mbox{і}\quad
\|h\|_{\alpha,\Gamma}:=(h,h)_{\alpha,\Gamma}^{1/2}.
\end{equation*}

Тут, звісно, $\mathcal{D}'(\Gamma)$~--- лінійний топологічний простір усіх розподілів на многовиді $\Gamma$, антидуальний до лінійного топологічного простору $\mathcal{D}(\Gamma)=C^{\infty}(\Gamma)$ усіх нескінченно диференційовних функцій на $\Gamma$, причому дуальність розглядається відносно скалярного добутку в гільбертовому просторі $L_{2}(\Gamma,dx)$ функцій, квадратично інтегровних на $\Gamma$ відносно $dx$. Вираз $(\chi_{j}h)\circ\pi_{j}$ позначає зображення розподілу $h$ у локальній карті~$\pi_{j}$. З~означення випливають неперервні вкладення $\mathcal{D}(\Gamma)\hookrightarrow H^\alpha(\Gamma)\hookrightarrow \mathcal{D}'(\Gamma)$.

\begin{theorem}
Простір $H^\alpha(\Gamma)$ гільбертів і сепарабельний. Він з точністю до
еквівалентності норм не залежить від зазначеного вибору атласу і
розбиття одиниці. Множина $C^{\infty}(\Gamma)$ щільна у цьому просторі.
\end{theorem}

Клас просторів $\{H^{\alpha}(\Gamma):\alpha\in\mathrm{RO}\}$
називаємо розширеною соболєвською шкалою на $\Gamma$. Він уведений і досліджений в \cite{MikhailetsMurach09Dop3} (див. також монографію \cite[п.~2.4.2]{MikhailetsMurach14}).

Якщо $\alpha(t)\equiv t^{s}$ для деякого $s\in\mathbb{R}$, то  $H^{\alpha}(\Gamma)$~---  гільбертів простір Соболєва $H^{(s)}(\Gamma)$ порядку $s$ на $\Gamma$. Зокрема, $H^{(0)}(\Gamma)=L_{2}(\Gamma,dx)$ з еквівалентністю норм.

\begin{theorem}
Теорема~$\ref{th-properties-H}$ і властивість $\eqref{21.6}$ зберігають силу для відповідних просторів на $\Gamma$. Крім того, вкладення $H^{\alpha_{1}}(\Gamma)\hookrightarrow H^{\alpha}(\Gamma)$ компактне тоді і тільки тоді, коли $\alpha(t)/\alpha_{1}(t)\to\nobreak0$ при $t\to\infty$. Тому є компактними вкладення в аналогах твердження~$(\mathrm{iii})$ і властивості $\eqref{21.6}$ для $\Gamma$.
\end{theorem}

\begin{theorem}
Теореми $\ref{th-int-Sobolev}$--$\ref{th-int-space}$ зберігають силу, якщо в них замінити $\Omega$ на $\Gamma$.
\end{theorem}

За означенням, $\|w\|_{\alpha,\mathbb{R}^{n}}=
\|\alpha((1-\Delta)^{1/2})w\|_{\mathbb{R}^{n}}$ для довільного $w\in C^{\infty}_{0}(\mathbb{R}^{n})$, де $\|\cdot\|_{\mathbb{R}^{n}}$~--- норма в $L_{2}(\mathbb{R}^{n})$. Тут, як звичайно, $\Delta$~--- оператор Лапласа в $\mathbb{R}^{n}$, а $\alpha((1-\Delta)^{1/2})$~--- борелева функція $\alpha(t^{1/2})$ аргументу $t\geq1$ від самоспряженого необмеженого оператора $1-\Delta$ в $L_{2}(\mathbb{R}^{n})$. Подібний результат правильний для простору $H^{\alpha}(\Gamma)$. Припустимо, що многовид $\Gamma$ наділено рімановою метрикою, узгодженою зі щільністю $dx$ (це завжди можна зробити). На такому многовиді означено диференціальний оператор Бельтрамі--Лапласа $\Delta_{\Gamma}$. Він є додатним самоспряженим необмеженим оператором в просторі $L_{2}(\Gamma,dx)$ з областю визначення $H^{(2)}(\Gamma)$.

\begin{theorem}
Для кожного $\alpha\in\mathrm{RO}$ норми $\|u\|_{\alpha,\Gamma}$ і
$\|\alpha((1-\Delta_{\Gamma})^{1/2})u\|_{\Gamma}$ еквівалентні на класі всіх функцій $u\in C^{\infty}(\Gamma)$; тут $\|\cdot\|_{\Gamma}$~--- норма у просторі $L_{2}(\Gamma,dx)$. Отже, простір $H^{\alpha}(\Gamma)$ є (з точністю до еквівалентності норм) поповненням лінійного простору $C^{\infty}(\Gamma)$ за другою нормою.
\end{theorem}

Ці властивості розширеної соболєвської шкали на $\Gamma$ наведені в
\cite[пп. 2,~4]{MikhailetsMurach09Dop3} і \cite[п.~2.4.2]{MikhailetsMurach14}).

Якщо $\Omega$~--- обмежена (відкрита) область в $\mathbb{R}^{n}$ класу $C^{\infty}$, де ціле $n\geq2$, то її межа $\partial\Omega$ є нескінченно гладким компактним многовидом без краю і вимірності $n-1$ з $C^{\infty}$-структурою, успадкованою з простору $\mathbb{R}^{n}$. Як звичайно, $\overline{\Omega}:=\Omega\cup\partial\Omega$~--- замикання області~$\Omega$. Між розширеними соболєвськими шкалами на $\Omega$ і $\partial\Omega$ існує такий зв'язок.

\begin{theorem}
Нехай $\alpha\in\mathrm{RO}$ і $\sigma_0(\alpha)>0$. Тоді відображення $R:u\mapsto u\!\upharpoonright\!\partial\Omega$, де $u\in C^\infty(\overline{\Omega})$, продовжується за неперервністю до обмеженого лінійного оператора сліду
\begin{equation*}%%\label{trace}
R:H^{\alpha\varrho^{1/2}}(\Omega)\rightarrow
H^{\alpha}(\partial\Omega),
\end{equation*}
де $\alpha\varrho^{1/2}$ позначає функцію $\alpha(t)t^{1/2}$ аргументу $t\geq1$. Цей оператор сюр'єктивний і виконується еквівалентність норм
\begin{equation*}
\|h\|_{\alpha,\partial\Omega}\asymp \inf\bigl\{\|u\|_{\alpha\varrho^{1/2},\Omega}:\,u\in H^{\alpha\varrho^{1/2}}(\Omega),\; Ru=h\bigr\},
\end{equation*}
де $h\in H^\alpha(\Gamma)$.
\end{theorem}

Цей результат випливає з соболєвського випадку, коли $\alpha(t)\equiv t^{s}$ і $s>0$, за допомогою інтерполяційної теореми~\ref{th-int-Sobolev} та її аналогу для просторів на $\partial\Omega$.

Розширена соболєвська шкала на нескінченно гладкому компактному многовиді з краєм уведена і досліджена в \cite[пп. 1,~2]{KasirenkoMurachChepurukhina19Dop3}. Її властивості подібні до властивостей розглянутих вище шкал. Уточнені соболєвські шкали на компактному многовиді $\Gamma$ (з краєм або без) уведені та досліджені в
\cite[п.~3]{MikhailetsMurach06UMJ3} і \cite[п.~3.3]{MikhailetsMurach08MFAT1}. Вони складаються з просторів
$H^{s,\varphi}(\Gamma):=H^{\alpha}(\Gamma)$, для яких показник гладкості $\alpha\in\mathrm{RO}$ набирає вигляду $\alpha(t)\equiv t^{s}\varphi(t)$, де $s\in\mathbb{R}$, а функція $\varphi$ є повільно змінною на нескінченності за Караматою. Уточнена соболєвська шкала на векторних розшаруваннях над нескінченно гладким компактним многовидом без краю уведена і досліджена в \cite[пп. 2, 4,~7]{Zinchenko17OpenMath}, відповідна розширена соболєвська шкала~--- в \cite[пп. 2,~4]{Zinchenko17Collection3}.

Якщо $\Gamma$~--- коло, то уточнена і розширена соболєвські шкали на $\Gamma$ тісно пов'язані \cite[п.~2.1.3, приклад~2.1]{MikhailetsMurach14} з просторами, уведеними О.~І.~Степанцем \cite[част.~I, розд.~3, п.~7.1]{Stepanets05} і застосованими у теорії апроксимації.

\section{Еліптичні системи в евклідовому просторі}\label{survey-sect6}

Нехай $n,p\in\mathbb{N}$. Розглядається система $p$ лінійних диференціальних рівнянь в евклідовому просторі $\mathbb{R}^{n}$:
\begin{equation}\label{sect6f1}
\sum _{k=1}^p A_{j,k}u_{k}=f_{j},\quad j=1,\ldots,p,
\end{equation}
де
$$
A_{j,k}:=A_{j,k}(x,D):=\sum_{|\mu|\leq r_{j,k}}a_{\mu}^{j,k}(x)D^\mu,\quad j,k=1,\ldots,p.
$$
Тут і надалі $\mu=(\mu_1,\ldots,\mu_n)$~--- мультиіндекс з невід'ємними цілими компонентами, $|\mu|$:= $\mu_1+\ldots+\mu_n$, $D_{j}:=i\partial/\partial x_{j}$,
$D^\mu:=D^\mu_x:=D_{1}^{\mu_1}\ldots D_{n}^{\mu_n}$, де $i$~--- уявна одиниця, а $x=(x_{1},\ldots,x_{n})\in\mathbb{R}^n$. Перетворення Фур'є
переводить диференціальний оператор $D^\mu$ в оператор множення на функцію $\xi^\mu:=\xi_{1}^{\mu_1}\ldots\xi_{n}^{\mu_n}$ аргументу
$\xi=(\xi_{1},\ldots,\xi_{n})\in\mathbb{R}^n$, дуального до просторової змінної~$x$. Припускаємо, що кожний коефіцієнт $a_{\mu}^{j,k}\in C^{\infty}_{\mathrm{b}}(\mathbb{R}^n)$, тобто він є нескінченно диференційовною функцію на $\mathbb{R}^n$ та його довільна частинна похідна обмежена на~$\mathbb{R}^n$.

Розв'язок системи \eqref{sect6f1} розуміємо у сенсі теорії розподілів; отже, частинні похідні $D^{\mu}u_{k}$ у цій системі є узагальненими. Запишемо систему у матричній формі $\mathcal{A}u=f$, де  $\mathcal{A}:=(A_{j,k})_{j,k=1}^{p}$, $u=\mathrm{col}\,(u_{1},\ldots,u_{p})$, $f=\mathrm{col}\,(f_{1},\ldots,f_{p})$.

Припускаємо, що система \eqref{sect6f1} є \textit{рівномірно еліптичною} в $\mathbb{R}^{n}$ за Дуглісом--Ніренбергом \cite[п. 3.2.b]{Agranovich94}, тобто існують набори цілих чисел $\{l_1,\ldots,l_p\}$ і $\{m_1,\ldots,m_p\}$ такі, що:
\begin{itemize}
  \item[(i)] $r_{j,k}\leq l_j+m_k$ для усіх $j,k\in\{1,\ldots,p\}$ (якщо $l_j+m_k<0$, то $A_{j,k}\equiv0$);
  \item[(ii)] для деякого числа $c>0$ виконується нерівність $|\det(A_{j, k}^{\circ}(x,\xi))_{j,k=1}^p|\geq c$ для довільних $x,\xi\in\mathbb{R}^{n}$ із $|\xi|=1$.
\end{itemize}
Тут
\begin{equation*}
A_{j, k}^{\circ}(x,\xi):=\sum_{|\mu|= l_j+m_k}a_{\mu}^{j,k}(x)\,\xi^\mu
\end{equation*}
--- головний символ диференціального оператора $A_{j,k}(x,D)$ у випадку
$r_{j,k}=l_j+m_k$, або $A_{j, k}^{\circ}(x,\xi)\equiv0$ у випадку $r_{j,k}<l_j+m_k$.

Якщо всі $l_{j}=0$, то система \eqref{sect6f1} називається рівномірно еліптичною за Петровським. Якщо, крім того, всі числа $m_k$ однакові, то вона називається рівномірно еліптичною у звичайному сенсі.

Матричний диференціальний оператор $\mathcal{A}$ є обмеженим на парі гільбертових просторів
\begin{equation}\label{operator-A-Rn}
\mathcal{A}:\bigoplus_{k=1}^pH^{\varphi\varrho^{m_k}}(\mathbb{R}^{n})\to
\bigoplus_{j=1}^p\,H^{\varphi\varrho^{-l_j}}(\mathbb{R}^{n})
\end{equation}
для кожного $\varphi\in\mathrm{RO}$. Тут і надалі, $\varphi\varrho^{s}$ позначає функцію $\varphi(t)t^{s}$ аргументу $t\geq1$, де $s\in\mathbb{R}$. Це позначення використовуємо для того, щоб не писати аргумент $t$ у показниках гладкості функціональних просторів. Звісно, $\varphi\varrho^{s}\in\mathrm{RO}$ і $\sigma_{r}(\varphi\varrho^{s})=\sigma_{r}(\varphi)+s$ для кожного $r\in\{0,1\}$.

\begin{theorem}\label{survey-th6.1}
Нехай $\varphi\in\mathrm{RO}$, а число $\sigma>0$. Тоді існує число $c=c(\varphi,\sigma)>0$ таке, що для довільних вектор-функцій
\begin{equation*}
u\in\bigoplus_{k=1}^{p}H^{\varphi\varrho^{m_k}}(\mathbb{R}^{n})
\quad\mbox{і}\quad
f\in\bigoplus_{j=1}^{p}H^{\varphi\varrho^{-l_j}}(\mathbb{R}^{n}),
\end{equation*}
які задовольняють рівняння $\mathcal{A}u=f$ в $\mathbb{R}^{n}$, виконується апріорна оцінка
\begin{equation*}
\sum_{k=1}^{p}\|u_{k}\|_{\varphi\varrho^{m_k},\mathbb{R}^{n}}\leq
c\biggl(\,\sum_{j=1}^{p}\|f_{j}\|_{\varphi\varrho^{-l_j},\mathbb{R}^{n}}+
\sum_{k=1}^{p}\|u_{k}\|_{\varphi\varrho^{m_k-\sigma},\mathbb{R}^{n}}\biggr).
\end{equation*}
\end{theorem}

Позначимо через $H^{-\infty}(\mathbb{R}^{n})$ об'єднання усіх просторів $H^{(s)}(\mathbb{R}^{n})$, де $s\in\mathbb{R}$ (воно збігається з об'єднанням усіх просторів $H^{\alpha}(\mathbb{R}^{n})$, де $\alpha\in\mathrm{RO}$). Нехай $V$~--- довільна відкрита непорожня підмножина простору $\mathbb{R}^{n}$. Позначимо через $H^{\alpha}_{\mathrm{int}}(V)$, де $\alpha\in\mathrm{RO}$, лінійний простір усіх розподілів $w\in H^{-\infty}(\mathbb{R}^{n})$ таких, що $\chi w\in H^{\alpha}(\mathbb{R}^{n})$ для кожної функції $\chi\in C^{\infty}_{\mathrm{b}}(\mathbb{R}^{n})$, яка задовольняє умови $\mathrm{supp}\,\chi\subset V$ і $\mathrm{dist}(\mathrm{supp}\,\chi,\partial
V)>0$. Як звичайно, $\mathrm{supp}\,\chi$~--- носій функції (або розподілу) $\chi$, $\partial V$~--- межа множини $V$, а $\mathrm{dist}(\cdot,\cdot)$~--- відстань. Якщо $V=\mathbb{R}^{n}$, то $H^{\alpha}_{\mathrm{int}}(V)=H^{\alpha}(\mathbb{R}^{n})$.

\begin{theorem}\label{th-sect6-regularity}
Нехай $\varphi\in\mathrm{RO}$. Припустимо, що вектор-функція $u\in(H^{-\infty}(\mathbb{R}^{n}))^p$ є розв'язком системи $\mathcal{A}u=f$ на відкритій множині $V\subset\mathbb{R}^{n}$ та
\begin{equation}\label{survey-th6.2-cond-f}
f_j\in H_{\mathrm{int}}^{\varphi\varrho^{-l_j}}(V)
\quad\mbox{для кожного}\quad j\in\{1,\ldots,p\}.
\end{equation}
Тоді $u_{k}\in
H_{\mathrm{int}}^{\varphi\varrho^{m_k}}(V)$ для кожного $k\in\{1,\ldots,p\}$.
\end{theorem}

Слід розрізняти внутрішню і локальну гладкість (чи регулярність) на $V$. Остання характеризується у термінах просторів $H^{\alpha}_{\mathrm{loc}}(V)$, кожний з яких складається з усіх розподілів $w\in H^{-\infty}(\mathbb{R}^{n})$ таких, що $\chi w\in H^{\alpha}(\mathbb{R}^{n})$ для довільної функції $\chi\in C^{\infty}_{0}(\mathbb{R}^{n})$, яка задовольняє умову $\mathrm{supp}\,\chi\subset V$. Якщо множина $V$ обмежена, то
$H^{\alpha}_{\mathrm{int}}(V)=H^{\alpha}_{\mathrm{loc}}(V)$. Якщо
$V$ не обмежена, то можливе строге включення
$H^{\alpha}_{\mathrm{int}}(V)\subsetneqq H^{\alpha}_{\mathrm{loc}}(V)$. Для локальної гладкості виконується аналог останньої теореми; у її формулюванні слід лише замінити $\mathrm{int}$ на $\mathrm{loc}$. При цьому замість умови (ii) робиться більш слабке припущення, що $\det(A_{j, k}^{\circ}(x,\xi))_{j,k=1}^p\neq0$ для довільних $x\in V$ і $\xi\in\mathbb{R}^{n}$ (умова еліптичності матричного диференціального оператора $\mathcal{A}$ в кожній точці множини $V$).

\begin{theorem}\label{survey-th6.3}
Нехай $k,\lambda\in\mathbb{Z}$, $1\leq k\leq p$, $\lambda\geq0$, а параметр $\varphi\in\mathrm{RO}$ задовольняє умову
\begin{equation}\label{survey-th6.3-cond-C}
\int_{1}^{\,\infty}\;t^{2\lambda+n-1-2m_k}\varphi^{-2}(t)\,dt<\infty.
\end{equation}
Припустимо, що вектор-функція $u\in(H^{-\infty}(\mathbb{R}^{n}))^p$ задовольняє умову теореми~$\ref{th-sect6-regularity}$. Тоді всі узагальнені частинні похідні компоненти $u_k$ розв'язку до порядку $\lambda$ включно неперервні на множині $V$, причому ці похідні обмежені на кожній множині $V_{0}\subset V$ такій, що  $\mathrm{dist}(V_{0},\partial V)>0$. Зокрема, якщо
$V=\mathbb{R}^{n}$, то $u_k\in C_{\mathrm{b}}^{\lambda}(\mathbb{R}^{n})$. \end{theorem}

\begin{remark}\label{survey-rem6.1}
У теоремі \ref{survey-th6.3} умова \eqref{survey-th6.3-cond-C} є точною, тобто ця умова еквівалентна такій властивості:
\begin{equation}
\bigl(u\in(H^{-\infty}(\mathbb{R}^{n}))^p,\;\mathcal{A}u=f,
\;f\,\mbox{задовольняє \eqref{survey-th6.2-cond-f}}\bigr)\;\Longrightarrow\;
u_k\in C^{\lambda}(V).
\end{equation}
\end{remark}

Якщо в умові \eqref{survey-th6.2-cond-f} замінити $\mathrm{int}$ на $\mathrm{loc}$ в позначенні простору, то з умов  теореми~\ref{survey-th6.3} випливатиме лише, що $u_k\in C^{\lambda}(V)$. При цьому достатньо припускати еліптичність оператора $\mathcal{A}$ в кожній точці множини $V$. За цього припущення виконується

\begin{corollary}\label{cor3.2}
Нехай усі $l_j=0$. Припустимо, що
вектор-функція $u\in H^{-\infty}(\mathbb{R}^{n})$ є розв'язком системи $\mathcal{A}u=f$ на
відкритій множині $V\subseteq\mathbb{R}^{n}$, де
$f\in(H_{\mathrm{loc}}^{\varphi}(V))^p$, а параметр
$\varphi\in\mathrm{RO}$ задовольняє умову
\begin{equation*}
\int_{1}^{\infty}\,t^{n-1}\,\varphi^{-2}(t)\,dt<\infty.
\end{equation*}
Тоді $u_{k}\in C^{m_k}(V)$ для кожного $k\in\{1,\ldots,m\}$, тобто розв'язок $u$ є класичним на~$V$.
\end{corollary}

Для класичного розв'язку $u$ вектор-функція $f=\mathcal{A}u$ обчислюється на $V$ за допомогою класичних частинних похідних, оскільки на кожну компоненту $u_{k}\in C^{m_k}(V)$ розв'язку діють диференціальні оператори $A_{j,k}$ порядки яких не перевищують $m_k$.

Аналоги цих теорем правильні для системи $\mathcal{A}^{+}v=g$, формально
спряженої до системи \eqref{sect6f1}, оскільки обидві вони рівномірно еліптичні в $\mathbb{R}^{n}$ за Дуглісом--Ніренбергом. Тут
$\mathcal{A}^{+}:=(A_{k,j}^{+}(x,D))_{j,k=1}^{p}$, де
\begin{equation*}
A_{k,j}^{+}(x,D)v_{k}(x):=\sum_{|\mu|\leq
r_{k,j}}D^{\mu}\bigl({\overline{a_{\mu}^{k,j}(x)}}v_{k}(x)\bigr).
\end{equation*}
Отже, $(\mathcal{A}^{+}u,v)_{\mathbb{R}^n}=(u,\mathcal{A}v)_{\mathbb{R}^n}$ для довільних $u,v\in(\mathcal{S}(\mathbb{R}^{n}))^p$, де  $(\cdot,\cdot)_{\mathbb{R}^n}$~--- скалярний добуток у гільбертовому просторі $(L_{2}(\mathbb{R}^{n}))^p$. Матричний диференціальний оператор  $\mathcal{A}^{+}$ є обмеженим на парі просторів
\begin{equation}\label{operator-A-adjoint-Rn}
\mathcal{A}^{+}:\bigoplus_{k=1}^p H^{\varphi\varrho^{l_k}}(\mathbb{R}^{n})\rightarrow
\bigoplus_{j=1}^p H^{\varphi\varrho^{-m_j}}(\mathbb{R}^{n})
\end{equation}
для кожного $\varphi\in\mathrm{RO}$.

За теоремою~\ref{th-sect6-regularity} ядра операторів \eqref{operator-A-Rn} і \eqref{operator-A-adjoint-Rn} збігаються з просторами
\begin{equation*}
\mathcal{N}:=\bigl\{u\in(H^{\infty}(\mathbb{R}^{n}))^p:\mathcal{A}u=0\bigr\}
\quad\mbox{і}\quad
\mathcal{N}^{+}:=
\bigl\{v\in(H^{\infty}(\mathbb{R}^{n}))^p:\mathcal{A}^{+}v=0\bigr\}
\end{equation*}
відповідно й тому не залежать від $\varphi$. Тут  $H^{\infty}(\mathbb{R}^{n})$ позначає перетин усіх просторів $H^{(s)}(\mathbb{R}^{n})$, де $s\in\mathbb{R}$; цей перетин лежить в $C^{\infty}_{\mathrm{b}}(\mathbb{R}^{n})$.

Система наступних двох умов є достатньою для фредгольмовості оператора \eqref{operator-A-Rn} (а також оператора \eqref{operator-A-adjoint-Rn}):
\begin{itemize}
  \item[(a)] $D^{\beta}a_{\mu}^{j,k}(x)\rightarrow0$ при $|x|\to\infty$ для кожного мультиіндексу $\beta$ з $|\beta|\geq 1$, довільних індексів
$j,k\in\{1,\ldots,p\}$ і мультиіндексу $\mu$ з $|\mu|\leq r_{j,k}$;
  \item[(b)] існують числа $c_{1}>0$ і $c_{2}\geq0$ такі, що $|\det(A_{j, k}(x,\xi))_{j,k=1}^p|\geq c_{1}\langle\xi\rangle^{q}$ для довільних векторів $x,\xi\in\mathbb{R}^{n}$, які задовольняють $|x|+|\xi|\geq c_{2}$; тут $q:=l_1+\ldots+l_p+m_1+\ldots+m_p$, а
\begin{equation*}
A_{j, k}(x,\xi):=\sum_{|\mu|\leq l_j+m_k}a_{\mu}^{j,k}(x)\,\xi^\mu
\end{equation*}
---~повний символ диференціального оператора $A_{j,k}(x,D)$.
\end{itemize}

\begin{theorem}\label{survey-th6.4}
Нехай виконуються умови $\mathrm{(a)}$ і $\mathrm{(b)}$. Тоді оператор \eqref{operator-A-Rn} фредгольмів для кожного $\varphi\in\mathrm{RO}$. Область значень цього оператора складається з усіх вектор-функцій
\begin{equation*}
f=\mathrm{col}\,(f_{1},\ldots,f_{p})\in
\bigoplus_{j=1}^{p}H^{\varphi\varrho^{-l_j}}(\mathbb{R}^{n})
\end{equation*}
таких, що $(f_1,v_1)_{\mathbb{R}^{n}}+\ldots+(f_p,v_p)_{\mathbb{R}^{n}}=0$ для кожної вектор-функції $v=\mathrm{col}\,(v_{1},\ldots,v_{p})\in\mathcal{N}^{+}$. Отже, індекс оператора \eqref{operator-A-Rn} дорівнює $\dim\mathcal{N}-\dim\mathcal{N}^{+}$ й не залежить від~$\varphi$.
\end{theorem}

Тут $(\cdot,\cdot)_{\mathbb{R}^n}$ позначає півторалінійну форму відносно якої взаємно антидуальні простори $\mathcal{S}'(\mathbb{R}^{n})$ і $\mathcal{S}(\mathbb{R}^{n})$ (вона є розширенням за неперервністю скалярного добутку в $L_{2}(\mathbb{R}^{n})$).

Зауважимо, що $\mathrm{(b)}\Rightarrow\mathrm{(ii)}$. Якщо припустити (a), то умова (b) випливатиме з нетеровості оператора \eqref{operator-A-Rn} у випадку, коли $\varphi(t)\equiv1$ \cite[теорема 4.2]{Rabier12}.

Ці теореми доведено в \cite[п.~5]{ZinchenkoMurach12UMJ11} (а зауваження~\ref{survey-rem6.1} обґрунтовується подібно до \cite[зауваження~2]{Zinchenko17OpenMath}). Теореми \ref{survey-th6.1}--\ref{survey-th6.3} у скалярному випадку $p=1$ встановлені в \cite[пп. 5,~6]{MikhailetsMurach08Collection1}, а у випадку, коли усі $l_{j}=0$~--- в \cite{Murach09UMJ3}; для (більш вузької) уточненої соболєвської шкали вони доведені в \cite[пп. 4--6]{Murach08UMB3}; у цих роботах досліджені рівномірно еліптичні псевдодиференціальні оператори (див. також монографію \cite[пп.~1.4, 5.1]{MikhailetsMurach14}). У~соболєвському випадку $\varphi(t)\equiv t^{s}$ теореми \ref{survey-th6.1} і \ref{th-sect6-regularity} (для $V=\mathbb{R}^{n}$) добре відомі (див., наприклад, огляд М.~С.~Аграновича \cite[пп. 1.8, 3.2b]{Agranovich94}), а теорема~\ref{survey-th6.4} випливає з результату П.~Дж.~Рабієра \cite[теорема~4.2]{Rabier12}.

Локальна регулярність розв'язків гіпоеліптичних диференціальних рівнянь у просторах $\mathcal{B}_{p,k}(\mathbb{R}^{n})$
досліджена у монографіях Л.~Хермандера \cite[теореми 4.1.5, 7.4.1]{Hermander63} і \cite[теореми 11.1.8, 13.4.1]{Hermander05v2}. У  праці \cite[п.~3.8]{Hermander63} розглянуто також випадок систем  диференціальних рівнянь.

Скалярні еліптичні диференціальні рівняння на~$\mathbb{R}^{n}$ зі спектральним параметром досліджені недавно М.~Файерманом \cite{Faierman20} в $L_{p}$-просторах Соболєва, для яких показником гладкості служить досить загальна функція $n$ частотних змінних. У цій роботі знайдено достатні умови, за яких виконується апріорна оцінка розв'язку рівняння і воно має єдиний розв'язок у вказаних просторах, а також достатні умови фредгольмовості  відповідного диференціального оператора на парах цих просторів \cite[теореми 3.4, 3.5, 4.1]{Faierman20}. Якщо $p=2$, то ці простори потрапляють у клас просторів $\mathcal{B}_{2,k}(\mathbb{R}^{n})$; якщо $p\neq2$, то вони відносяться до просторів, уведених і досліджених Л.~Р.~Волевичем і Б.~П.~Панеяхом \cite[\S~13]{VolevichPaneah65}.

У~монографії Ф.~Нікола і Л.~Родіно \cite{NicolaRodino10} викладено числення широких класів гіпоеліптичних псевдодиференціальних операторів в $\mathbb{R}^{n}$ у функціональних нормованих просторах змінної узагальненої гладкості (яка залежить не лише від частотних змінних, а також від просторових координат). Такі простори природно пов'язані із символом оператора і є широким узагальненням просторів Соболєва і Хермандера. За деяких умов глобальної еліптичності доведено  теореми про регулярність розв'язків псевдодиференціальних рівнянь та встановлено фредгольмовість відповідних операторів. Досліджено також спектральні властивості цих операторів (найголовніше, асимптотичну поведінку їх власних чисел).

У статті О.~Мілатовича \cite{Milatovic24} розглянуто
скалярні псевдодиференціальні оператори в уведеній ним розширеній соболєвській шкалі на решітці $\mathbb{Z}^{n}$. Доведено, що за деякої умови еліптичності оператора він є фредгольмовим на парі гільбертових просторів $H^{\varphi}(\mathbb{Z}^{n})$ і $H^{\varphi\varrho^{-m}}(\mathbb{Z}^{n})$, де $\varphi\in\mathrm{RO}$, a $m$~--- порядок оператора, причому індекс оператора не залежить від $\varphi$ \cite[теорема~2.6]{Milatovic24}.

\section{Еліптичні системи на замкненому многовиді}\label{survey-sect7}

Нехай $\Gamma$~--- многовид вимірності $n\in\mathbb{N}$, зазначений в розд.~\ref{survey-sect5}. На $\Gamma$ розглядається система $p\in\mathbb{N}$
лінійних диференціальних рівнянь вигляду \eqref{sect6f1} тобто $\mathcal{A}u=f$, де $\mathcal{A}:=(A_{j,k})_{j,k=1}^p$, а кожне $A_{j,k}$~--- скалярний лінійний диференціальний оператор на $\Gamma$ деякого порядку $r_{j,k}$ з коефіцієнтами класу $C^\infty(\Gamma)$. Розв'язок $u$ цієї системи розуміємо у сенсі теорії розподілів на $\Gamma$.

Припускаємо, що система $\mathcal{A}u=f$ є \textit{еліптичною} на $\Gamma$ за Дуглісом--Ніренбергом \cite[п.~3.2.b]{Agranovich94}, тобто існують набори цілих чисел $\{l_1,\ldots,l_p\}$ і $\{m_1,\ldots,m_p\}$ такі, що виконується умова~(i) з розд.~\ref{survey-sect6} та нерівність $\det(A_{j, k}^{\circ}(x,\xi))_{j,k=1}^p \neq 0$ для довільної точки
$x\in\Gamma$ і ковектора $\xi\in T_x^*\Gamma\setminus\nobreak\{0\}$. Тут $A_{j, k}^{\circ}(x,\xi)$~--- головний символ диференціального оператора $A_{j,k}$ у випадку $r_{j,k}=l_j+m_k$, інакше $A_{j,k}^{\circ}(x,\xi)\equiv0$. Як звичайно, $T_x^*\Gamma$~--- кодотичний простір до многовиду $\Gamma$ в точці~$x$.

Позначимо через $\mathcal{N}$ простір усіх вектор-функцій $u\in(C^\infty(\Gamma))^p$ таких, що $\mathcal{A}u=0$ на $\Gamma$. Нехай $\mathcal{A}^{+}$~--- матричний диференціальний оператор, формально спряжений до $\mathcal{A}$ відносно скалярного добутку в гільбертовому просторі $(L_{2}(\Gamma,dx))^{p}$. Позначимо через $\mathcal{N}^{+}$ простір усіх вектор-функцій $v\in(C^\infty(\Gamma))^p$ таких, що $\mathcal{A}^{+}v=0$ на $\Gamma$. Оскільки системи $\mathcal{A}u=f$ і $\mathcal{A}^{+}v=g$ еліптичні на $\Gamma$, то простори $\mathcal{N}$ і $\mathcal{N}^{+}$ скінченновимірні.

\begin{theorem}\label{survey-th7.1}
Оператор $\mathcal{A}$ є обмеженим і фредгольмовим на парі гільбертових просторів
\begin{equation}\label{operator-A-Gamma}
\mathcal{A}:\bigoplus_{k=1}^p\,H^{\varphi\varrho^{m_k}}(\Gamma)\rightarrow
\bigoplus_{j=1}^p\,H^{\varphi\varrho^{-l_j}}(\Gamma)
\end{equation}
для кожного $\varphi\in\mathrm{RO}$. Його ядро дорівнює $\mathcal{N}$, а область значень складається з усіх вектор-функцій
\begin{equation*}
f=\mathrm{col}\,(f_{1},\ldots,f_{p})\in
\bigoplus_{j=1}^{p}H^{\varphi\varrho^{-l_j}}(\Gamma)
\end{equation*}
таких, що $(f_1,v_1)_{\Gamma}+\ldots+(f_p,v_p)_{\Gamma}=0$ для кожної вектор-функції $v=\mathrm{col}\,(v_{1},\ldots,v_{p})\in\mathcal{N}^{+}$. Отже, індекс оператора \eqref{operator-A-Gamma} дорівнює $\dim\mathcal{N}-\dim\mathcal{N}^{+}$ й не залежить від~$\varphi$.
\end{theorem}

Тут $(\cdot,\cdot)_{\Gamma}$ позначає півторалінійну форму, відносно якої взаємно антидуальні простори $\mathcal{D}'(\Gamma)$ і $\mathcal{D}(\Gamma)$ (вона є розширенням за неперервністю скалярного добутку в $L_{2}(\Gamma,dx)$). Якщо $p=1$, то
індекс оператора \eqref{operator-A-Gamma} дорівнює нулю.

Нехай $\Gamma_{0}$~---довільна відкрита непорожня підмножина многовиду $\Gamma$. Позначимо через $H_{\mathrm{loc}}^{\alpha}(\Gamma_{0})$, де $\alpha\in\mathrm{RO}$, лінійний простір усіх розподілів $h\in\mathcal{D}'(\Gamma)$ таких, що
$\chi h\in H^{\alpha}(\Gamma)$ для кожної функції $\chi\in C^\infty(\Gamma)$, яка задовольняє умову $\mathrm{supp}\,\chi\subset\Gamma_{0}$. Якщо $\Gamma_{0}=\Gamma$, то
$H_{\mathrm{loc}}^{\alpha}(\Gamma_{0})=H^{\alpha}(\Gamma)$.

\begin{theorem}\label{survey-th7.2}
Нехай $\varphi\in\mathrm{RO}$. Припустимо, що вектор-функція
$u\in(\mathcal{D}'(\Gamma))^p$ є розв'язком системи $\mathcal{A}u=f$ на відкритій множині $\Gamma_{0}\subset\Gamma$ та $f_{j}\in
H_{\mathrm{loc}}^{\varphi\varrho^{-l_j}}(\Gamma_{0})$ для кожного $j\in\{1,\ldots,p\}$. Тоді $u_{k}\in
H_{\mathrm{loc}}^{\varphi\varrho^{m_k}}(\Gamma_{0})$ для кожного $k\in\{1,\ldots,p\}$.
\end{theorem}

Цю теорему доповнюють такі результати:

\begin{theorem}\label{survey-th7.3}
Нехай $\varphi\in\mathrm{RO}$. Припустимо, що вектор-функція
$u\in(\mathcal{D}'(\Gamma))^p$ задовольняє умову теореми~$\ref{survey-th7.2}$. Довільно виберемо число $r>0$ і функції $\chi,\chi_{1}\in C^\infty(\Gamma)$ такі, що $\mathrm{supp}\,\chi\subset \mathrm{supp}\,\chi_{1}\subset\Gamma_{0}$ і $\chi_{1}=1$ в околі $\mathrm{supp}\,\chi$. Тоді існує число $c=c(\varphi,r,\chi,\chi_{1})>0$ таке, що
\begin{equation*}
\sum_{k=1}^{p}\|\chi u_{k}\|_{\varphi\varrho^{m_k},\Gamma}\leq
c\biggl(\,\sum_{j=1}^{p}\|\chi_{1}f_{j}\|_{\varphi\varrho^{-l_j},\Gamma}+
\sum_{k=1}^{p}\|\chi_{1} u_{k}\|_{\varphi\varrho^{m_k-r},\Gamma}\biggr),
\end{equation*}
причому $c$ не залежить від $u$ і $f=\mathcal{A}u$. Якщо $0<r\leq1$, то можна узяти $\chi f_{j}$ замість $\chi_{1}f_{j}$ в оцінці.
\end{theorem}

\begin{theorem}\label{survey-th7.4}
Нехай $k,\lambda\in\mathbb{Z}$, $1\leq k\leq p$, $\lambda\geq0$, а параметр $\varphi\in\mathrm{RO}$ задовольняє умову~\eqref{survey-th6.3-cond-C}. Припустимо, що вектор-функція $u\in(\mathcal{D}'(\Gamma))^p$ задовольняє умову теореми~$\ref{survey-th7.2}$. Тоді $u_k\in C^\lambda(\Gamma_{0})$.
\end{theorem}

У теоремі~\ref{survey-th7.4} умова \eqref{survey-th6.3-cond-C} є точною у сенсі, аналогічному до зауваження~\ref{survey-rem6.1}. Якщо усі $l_{j}=0$, то з цієї теореми випливає достатня умова класичності розв'язку $u$, аналогічна наслідку~\ref{cor3.2}.

Ці теореми встановлені в \cite{Zinchenko13Dop3, Zinchenko14Collection2}; випадки одного рівняння і систем, еліптичних за Петровським, розглянуті в
\cite[п.~3]{MikhailetsMurach09Dop3} і \cite[пп. 3,~5]{ZinchenkoMurach14JMS5} відповідно. В уточненій соболєвській шкалі на $\Gamma$ еліптичні рівняння і системи (диференціальні та псевдодиференціальні) досліджені в \cite{MikhailetsMurach08BPAS3,
Murach07UMJ6, Murach08MFAT2, MikhailetsMurach06Dop10, Murach07Dop5} (див. також монографію \cite[пп. 2.2.2, 2.2.3, 5.2.2, 5.2.3]{MikhailetsMurach14}). У~цих працях версії теореми~\ref{survey-th7.3} встановлені у глобальному випадку, коли $V=\Gamma$. Загальний випадок досліджується подібно до
\cite[доведення теореми~8]{Zinchenko17OpenMath}. У цій праці \cite[пп. 5,~8]{Zinchenko17OpenMath} вивчено скалярні еліптичні псевдодиференціальні оператори в уточненій соболєвській шкалі на векторних розшаруваннях. У~соболєвському випадку $\varphi(t)\equiv t^{s}$ теореми \ref{survey-th7.1}--\ref{survey-th7.3} добре відомі (останні дві з них
принаймні для $\Gamma_0=\Gamma$; див., наприклад, огляд \cite[пп. 2.2, 2.3, 3.2b]{Agranovich94}).

Зауважимо, що уточнена і розширена соболєвські шкали на замкненому многовиді $\Gamma$ мають важливі застосування до дослідження умов різних типів збіжності розвинень за власними функціями нормальних (зокрема, самоспряжених) еліптичних диференціальних та псевдодиференціальних операторів на $\Gamma$. У термінах цих шкал вдається отримати точні достатні умови збіжності майже скрізь, безумовної збіжності майже скрізь та рівномірної збіжності вказаних розвинень і отримати тонкі оцінки швидкості рівномірної збіжності та збіжності у середньому разом з похідними вказаного порядку (див. монографію \cite[пп. 2.3.2, 2.3.3]{MikhailetsMurach14} і новітню статтю \cite{MikhailetsMurach24ProcRSEdinburgh}). Звісно, ці результати застосовні до кратних рядів Фур'є (випадок, коли $\Gamma$~--- багатовимірний тор), та до розвинень за власними функціями звичайних диференціальних операторів з періодичними коефіцієнтами ($\Gamma$~--- коло). В останньому випадку вказані шкали були використані у спектральній теорії операторів Хіла--Шрьодінгера \cite{MikhailetsMolyboga09, MikhailetsMolyboga11, MikhailetsMolyboga12}.

\section{Еліптичні системи з параметром}\label{survey-sect8}

Як і раніше, $\Gamma$~--- многовид вимірності $n\in\mathbb{N}$, зазначений у розд.~\ref{survey-sect5}. На $\Gamma$ розглядається система $p\in\mathbb{N}$ лінійних диференціальних рівнянь вигляду
\begin{equation}\label{survey-5.25}
\sum _{k=1}^p A_{j,k}(\lambda)u_{k}=f_{j},\quad j=1,\ldots,p.
\end{equation}
Тут $\lambda$~--- комплексний параметр, а кожне
\begin{equation*}
A_{j,k}(\lambda):=\sum_{r=0}^{m_{k}}\lambda^{r}A_{j,k}^{(r)},
\end{equation*}
де цілі числа $m_{1},\ldots,m_{p}\geq0$ фіксовані, а всі $A_{j,k}^{(r)}$~--- скалярні лінійні диференціальні оператори на $\Gamma$ з коефіцієнтами класу $C^\infty(\Gamma)$ і порядку $\mathrm{ord}\,A_{j,k}^{(r)}\leq m_{k}-r$. Запишемо цю систему в матричній формі $\mathcal{A}(\lambda)u=f$, де
$\mathcal{A}(\lambda)=(A_{j,k}(\lambda))_{j,k=1}^p$.

Пов'яжемо з оператором $\mathcal{A}(\lambda)$ матрицю
\begin{equation*}
\mathcal{A}^{(0)}(x,\xi,\lambda):= \biggl(\:\sum_{r=0}^{m_{k}}\,
\lambda^{r}A_{j,k}^{r,0}(x,\xi)\biggr)_{j,k=1}^{p},
\end{equation*}
яка залежить від $x\in\Gamma$, $\xi\in T_x^*\Gamma$ і  $\lambda\in\mathbb{C}$. Тут $A_{j,k}^{r,0}(x,\xi)$~--- головний символ
диференціального оператора $A_{j,k}^{(r)}$ у випадку, коли
$\mathrm{ord}\,A_{j,k}^{(r)}=m_{k}-r$, інакше $A_{j,k}^{r,0}(x,\xi)\equiv0$.

Нехай $K$~--- деякий фіксований замкнений кут на комплексній площині з вершиною у початку координат (він може вироджуватися у промінь). Припускаємо, що система \eqref{survey-5.25} є \textit{еліптичною з параметром} в куті $K$ на многовиді $\Gamma$, тобто $\det \mathcal{A}^{(0)}(x,\xi,\lambda)\neq0$ для довільних $x\in\Gamma$, $\xi\in
T_{x}^{\ast}\Gamma$ і $\lambda\in K$ таких, що $|\xi|+|\lambda|\neq0$.

Оскільки, для кожного фіксованого $\lambda\in\mathbb{C}$ ця система є еліптичною за Петровським на $\Gamma$, то $\mathcal{A}(\lambda)$ є фредгольмовим обмеженим оператором на парі просторів
\begin{equation}\label{survey-5.29}
\mathcal{A}(\lambda):\bigoplus_{k=1}^p H^{\varphi\varrho^{m_k}}(\Gamma)\rightarrow
(H^{\varphi}(\Gamma))^{p}
\end{equation}
для будь-якого $\varphi\in\mathrm{RO}$.

\begin{theorem}\label{survey-th8.1}
Існує число  $\lambda_{0}>0$ таке, що для довільних $\lambda\in K$ з $|\lambda|\geq\lambda_{0}$ і $\varphi\in\mathrm{RO}$ оператор \eqref{survey-5.29} є ізоморфізмом. Більше того, для кожного $\varphi\in\mathrm{RO}$ існує число $c=c(\varphi)\geq\nobreak1$ таке, що
\begin{equation*}
c^{-1}\,\sum_{j=1}^{p}\,\|f_{j}\|_{\varphi,\Gamma}\leq
\sum_{k=1}^{p}\,\bigl(\,\|u_{k}\|_{\varphi\varrho^{m_{k}},\Gamma}+
|\lambda|^{m_{k}}\|u_{k}\|_{\varphi,\Gamma}\,\bigr)\leq
c\,\sum_{j=1}^{p}\,\|f_{j}\|_{\varphi,\Gamma}
\end{equation*}
для кожного $\lambda\in K$ з $|\lambda|\geq\lambda_{0}$ і довільних вектор функцій
\begin{equation*}
u=\mathrm{col}(u_{1},\ldots,u_{p})\in
\bigoplus_{k=1}^p H^{\varphi\varrho^{m_k}}(\Gamma)
\quad\mbox{і}\quad
f=\mathrm{col}(f_{1},\ldots,f_{p})\in(H^{\varphi}(\Gamma))^{p},
\end{equation*}
які задовольняють рівняння $A(\lambda)u=f$. Число $c$ не залежить від $\lambda$, $u$ і $f$.
\end{theorem}

Звідси випливає, що індекс оператора \eqref{survey-5.29} дорівнює нулю для будь-яких $\lambda\in\mathbb{C}$ і $\varphi\in\mathrm{RO}$.

Ця теорема доведена в \cite[п.~4]{ZinchenkoMurach12Collection2}; скалярний випадок досліджено в \cite{MurachZinchenko08MFAT1}. В уточненій соболєвській шкалі різні класи еліптичних рівнянь і систем з параметром досліджені в \cite{Murach07UMJ6, Murach08MFAT2, Murach07Dop5} (відповідні результати викладені також у монографії \cite[пп. 2.2.4, 5.2.4]{MikhailetsMurach14}). У~соболєвському випадку $\varphi(t)\equiv t^{s}$ версії теореми~\ref{survey-th8.1} відомі для різних класів еліптичних рівнянь з параметром \cite[пп. 4.1,~4.3]{Agranovich94}.

\section{Еліптичні крайові задачі}\label{survey-sect9}

Надалі ціле $n\geq2$, $\Omega$~--- обмежена (відкрита) область в $\mathbb{R}^{n}$ класу $C^{\infty}$, $\Gamma$~--- межа цієї області ($\Omega$ лежить локально по один бік щодо $\Gamma$), $\nu(x)$~--- орт внутрішньої нормалі до $\Gamma$ у точці $x\in\Gamma$, а
$(\cdot,\cdot)_\Omega$ і $(\cdot,\cdot)_\Gamma$~--- скалярні добутки відповідно у просторах $L_2(\Omega)$ і $L_2(\Gamma,dx)$ та розширення за неперервністю цих скалярних добутків ($dx$~--- елемент площі на $\Gamma$).

Нехай цілі числа $q\geq1$ і $m_{1},\ldots,m_{q}$ такі, що $0\leq m_{j}\leq2q-1$ для кожного $j\in\{1,\ldots,q\}$. Покладемо $m:=\max\{m_{1},\dots,m_{q}\}$.

В області $\Omega$ розглядається лінійна крайова задача
\begin{gather}\label{survey-3f1}
Au(x)\equiv\sum_{|\mu|\leq{2q}}a_{\mu}(x)D^{\mu}u(x)=f(x),\quad x\in\Omega,\\
B_{j}u(x)\equiv \sum_{|\mu|\leq
m_{j}}b_{j,\mu}(x)D^{\mu}u(x)=g_{j}(x),\quad x\in\Gamma,\quad
j=1,\ldots,q. \label{survey-3f2}
\end{gather}
Тут усі $a_\mu\in C^{\infty}(\overline{\Omega})$ і $b_{j,\mu}\in C^\infty(\Gamma)$. Покладемо $B:=(B_{1},\ldots,B_{q})$.

Позначимо через $A^{\circ}(x,\xi)$ і $B_j^{\circ}(x,\xi)$ головні символи диференціальних операторів $A=A(x,D)$, де $x\in\overline{\Omega}$, і $B_j=B_j(x,D)$, де $x\in\Gamma$. Ці символи є однорідними поліномами аргументу $\xi\in\mathbb{C}^{n}$ степені $2q$ і $m_{j}$ відповідно.

Припускаємо, шо крайова задача  \eqref{survey-3f1}, \eqref{survey-3f2} є  \textit{регулярною еліптичною} в області $\Omega$, тобто задовольняє такі умови:
\begin{itemize}
\item [(i)] Диференціальний вираз $A$ є еліптичним на $\overline{\Omega}$, тобто $A^{\circ}(x,\xi)\neq0$ для будь-яких $x\in\Gamma$ і $\xi\in\mathbb{R}^{n}\setminus\{0\}$.
\item [(ii)] Вираз $A$ є правильно еліптичним на
$\Gamma$, тобто для кожної точки $x\in\Gamma$ і будь-якого ненульового вектора $\tau\in\mathbb{R}^{n}$, дотичного до $\Gamma$ у точці
$x\in\Gamma$, многочлен $A^{\circ}(x,\tau+\zeta\nu(x))$ аргументу $\zeta\in\mathbb{C}$ має $q$ коренів $\zeta^{+}_{j}(x,\tau)$, $j=1,\ldots,q$, з додатною уявною частиною і стільки ж коренів з від'ємною уявною частиною (з урахуванням їх кратності).
\item [(iii)] Система граничних диференціальних виразів $\{B_{1},\ldots,B_{q}\}$ задовольняє умову Шапіро--Лопатинського щодо $A$ на $\Gamma$, тобто для довільних $x$ і $\tau$, вказаних в умові~(ii), система многочленів $B_{j}^{\circ}(x,\tau+\zeta\nu(x))$, $j=1,\ldots,q$, аргументу $\zeta\in\mathbb{C}$ лінійно незалежна за модулем многочлена
    $\prod_{j=1}^{q}(\zeta-\zeta^{+}_{j}(x,\tau))$.
\item [(iv)] Система $\{B_{1},\ldots,B_{q}\}$ є нормальною на $\Gamma$ за Ароншайном, Мільграмом і Шехтером, тобто порядки $m_1,\ldots,m_q$ цих виразів парами різні і $B^{(0)}_j(x,\nu(x))\neq0$ для будь-яких $x\in\Gamma$ і  $j\in\{1,\ldots,q\}$.
\end{itemize}

Якщо $n\geq3$ або усі $a_\mu$ дійсні (при $n=2$), то $\mathrm{(i)}\Rightarrow\mathrm{(ii)}$.

Для опису області значень оператора, породженого задачею \eqref{survey-3f1}, \eqref{survey-3f2}, знадобиться формальна спряжена крайова задача
\begin{gather}\label{survey-deq373}
A^{+}v(x)\equiv\sum_{|\mu|\leq2q}D^{\mu}
\bigl(\,\overline{a_{\mu}(x)}\,v(x)\bigr)=
\omega(x),\quad x\in\Omega,\\
B^{+}_{j}v(x)=h_{j}(x),\quad x\in\Gamma,\quad j=1,\ldots,q. \label{survey-deq374}
\end{gather}
Спряження розглядається відносно формули Гріна
\begin{equation*}
(Au,v)_{\Omega} + \sum^{q}_{j=1}(B_{j}u,C^{+}_{j}v)_{\Gamma} =
(u,A^{+}v)_{\Omega} + \sum_{j=1}^{q}(C_{j}u,B^{+}_{j}v)_{\Gamma},
\end{equation*}
де $u,v\in C^{\infty}(\overline{\Omega})$. Тут $\{B^+_j\}$, $\{C_j\}$ і $\{C^+_j\}$~--- деякі нормальні системи граничних лінійних диференціальних виразів з коефіцієнтами класу $C^\infty(\Gamma)$. Порядки цих виразів задовольняють умову
$$
\mathrm{ord}\,B_j+\mathrm{ord}\,C^+_j=
\mathrm{ord}\,C_j+\mathrm{ord}\,B^+_j=2q-1.
$$

Позначимо через $N$ простір усіх розв'язків $u\in C^{\infty}(\overline{\Omega})$ задачі \eqref{survey-3f1}, \eqref{survey-3f2} з нульовими правими частинами. Аналогічно, $N^{+}$ позначає простір усіх розв'язків $v\in C^{\infty}(\overline{\Omega})$ задачі \eqref{survey-deq373}, \eqref{survey-deq374} з нульовими правими частинами. Ці простори скінченновимірні, оскільки обидві крайові задач еліптичні в~$\Omega$. Простір $N^{+}$ не залежить від вибору системи $\{B^{+}_{1},\ldots, B^{+}_{q}\}$, підпорядкованої формулі Гріна.

\begin{theorem}\label{survey-th9.1}
Нехай $\varphi\in\mathrm{RO}$ і $\sigma_0(\varphi)>m+1/2$. Тоді відображення $u\mapsto(Au,Bu)$, де $u\in C^\infty(\overline{\Omega})$,
продовжується єдиним чином (за неперервністю) до обмеженого лінійного оператора
\begin{equation}\label{survey-deq31}
(A,B):H^{\varphi}(\Omega)\rightarrow H^{\varphi\varrho^{-2q}}(\Omega)\oplus
\bigoplus_{j=1}^{q}H^{\varphi\varrho^{-m_j-1/2}}(\Gamma)=:
\mathcal{H}^{\varphi\varrho^{-2q}}(\Omega,\Gamma).
\end{equation}
Цей оператор фредгольмів. Його ядро збігається з $N$, а область
значень складається з усіх векторів $(f,g_1,\ldots,g_q)\in\mathcal{H}^{\varphi\varrho^{-2q}}(\Omega,\Gamma)$
таких, що
\begin{equation}\label{survey-range}
(f,v)_\Omega+\sum_{j=1}^{q}(g_j,C^{+}_{j}v)_{\Gamma}=0\quad\mbox{для довільного}\quad v\in N^{+}.
\end{equation}
Індекс оператора \eqref{survey-deq31} дорівнює $\dim N-\dim N^{+}$ й не залежить від~$\varphi$.
\end{theorem}

Позначимо через $H^{s+}(\Omega)$, де $s\in\mathbb{R}$, об'єднання усіх просторів $H^{(l)}(\Omega)$ з $l>s$ (воно збігається з об'єднанням усіх просторів $H^{\alpha}(\Omega)$, де $\alpha\in\mathrm{RO}$ і $\sigma_{0}(\alpha)>s$).

Розподіл $u\in H^{m+1/2+}(\Omega)$ називаємо
(сильним) \textit{узагальненим розв'язком} крайової задачі \eqref{survey-3f1}, \eqref{survey-3f2} з правою частиною
\begin{equation*}
(f,g):=(f,g_1,\ldots,g_q)\in
\mathcal{D}'(\Omega)\times(\mathcal{D}'(\Gamma))^q,
\end{equation*}
якщо $(A,B)u=(f,g)$ для оператора \eqref{survey-deq31} з деяким~$\alpha$ вказаним у теоремі~\ref{survey-th9.1}. Як звичайно, $\mathcal{D}'(\Omega)$ і $\mathcal{D}'(\Gamma)$~--- лінійні топологічні простори усіх розподілів, заданих відповідно в $\Omega$ і на $\Gamma$. Цей розв'язок має такі властивості.

Нехай $V$ --- довільна відкрита підмножина евклідового простору $\mathbb{R}^{n}$ така, що
$\Omega_0:=\Omega\cap V\neq\varnothing$. Покладемо
$\Gamma_{0}:=\Gamma\cap V$ (можливий випадок, коли
$\Gamma_{0}=\varnothing$). Позначимо через $H_{\mathrm{loc}}^{\alpha}(\Omega_{0},\Gamma_{0})$, де $\alpha\in\mathrm{RO}$, лінійний простір усіх розподілів $u\in
D'(\Omega)$ таких, що $\chi u\in H^{\alpha}(\Omega)$ для кожної функції $\chi\in C^{\infty}(\overline{\Omega})$, яка задовольняє умову $\mathrm{supp}\,\chi\subset\Omega_0\cup\Gamma_{0}$. Простір $H_{\mathrm{loc}}^{\alpha}(\Gamma_{0})$ означено в п.~\ref{survey-sect7}. Якщо $\Omega_0:=\Omega$ і $\Gamma_{0}:=\Gamma$, то $H_{\mathrm{loc}}^{\alpha}(\Omega_{0},\Gamma_{0})=H^{\alpha}(\Omega)$ і
$H_{\mathrm{loc}}^{\alpha}(\Gamma_{0})=H^{\alpha}(\Gamma)$.

\begin{theorem}\label{survey-th9.2}
Припустимо, що розподіл $u\in H^{m+1/2+}(\Omega)$ є узагальненим
розв'язком крайової задачі \eqref{survey-3f1}, \eqref{survey-3f2}, праві
частини якої задовольняють умови $f\in H^{\varphi\varrho^{-2q}}_{\mathrm{loc}}(\Omega_0,\Gamma_0)$ і $g_j\in H^{\varphi\varrho^{-m_j-1/2}}_{\mathrm{loc}}(\Gamma_0)$ для кожного номера $j\in\{1,\ldots,q\}$ і деякого параметра $\varphi\in\mathrm{RO}$ з $\sigma_0(\varphi)>m+1/2$. Тоді $u\in
H^{\varphi}_{\mathrm{loc}}(\Omega_0,\Gamma_0)$.
\end{theorem}

Якщо $\Gamma_{0}=\varnothing$, ця теорема стосується локальної регулярності розв'язку в околах внутрішніх точок області~$\Omega$.

\begin{theorem}\label{survey-th9.3}
Припустимо, що розподіл $u\in H^{m+1/2+}(\Omega)$ задовольняє умову теореми~$\ref{survey-th9.2}$. Виберемо довільно функції $\chi,\chi_{1}\in C^\infty(\overline{\Omega})$ такі, що
\begin{equation}\label{cutoff-functions}
\mathrm{supp}\,\chi\subset \mathrm{supp}\,\chi_{1}\subset\Omega_0\cup\Gamma_{0}
\quad\mbox{і}\quad
\chi_{1}=1\;\,\mbox{в околі}\;\,\mathrm{supp}\,\chi.
\end{equation}
Тоді існує число $c=c(\varphi,\chi,\chi_{1})>0$ таке, що
\begin{equation*}
\|\chi u\|_{\varphi,\Omega}\leq
c\biggl(\|\chi f\|_{\varphi\varrho^{-2q},\Omega}+
\sum_{j=1}^{q}\|\chi g_{j}\|_{\varphi\varrho^{-m_j-1/2},\Gamma}+
\|\chi_{1}u\|_{\varphi\varrho^{-1},\Omega}
\biggr),
\end{equation*}
причому $c$ не залежить від $u$ і $(f,g)=Au$.
\end{theorem}

\begin{theorem}\label{survey-th9.4}
Нехай ціле $\lambda\geq0$. Припустимо, що розподіл $u\in H^{m+1/2+}(\Omega)$ задовольняє умову теореми~$\ref{survey-th9.2}$, причому
\begin{equation}\label{survey-th9.4-cond}
\int_1^{\infty}t^{2\lambda+n-1}\varphi^{-2}(t)dt<\infty.
\end{equation}
Тоді $u\in C^{\lambda}(\Omega_{0}\cup\Gamma_{0})$.
\end{theorem}

Умова \eqref{survey-th9.4-cond} є точною у цій теоремі на класі зазначених розв'язків.

\begin{corollary}
Припустимо, що розподіл $u\in H^{m+1/2+}(\Omega)$ є узагальненим
розв'язком крайової задачі \eqref{survey-3f1}, \eqref{survey-3f2}, де
$f\in H^{\eta}_{\mathrm{loc}}(\Omega,\varnothing)\cap H^{\varphi\varrho^{-2q}}(\Omega)$ і $g_j\in H^{\varphi\varrho^{-m_j-1/2}}(\Gamma)$ для кожного номера $j\in\{1,\ldots,q\}$ та деяких параметрів $\eta,\varphi\in\mathrm{RO}$, які задовольняють умови $\sigma_0(\eta)>m+1/2$, $\sigma_0(\varphi)>m+1/2$ і
\begin{equation*}
\int_1^{\infty}t^{n-1}\eta^{-2}(t)dt<\infty,\quad
\int_1^{\infty}t^{2m+n-1}\varphi^{-2}(t)dt<\infty.
\end{equation*}
Тоді $u$~--- класичний розв'язок, тобто $u\in C^{2q}(\Omega)\cap C^{m}(\overline{\Omega})$.
\end{corollary}

Для класичного розв'язку $u$ крайової задачі \eqref{survey-3f1}, \eqref{survey-3f2} її ліві частини обчислюються за допомогою класичних частинних похідних та є неперервними функціями на $\Omega$ і $\Gamma$ відповідно.

Теореми \ref{survey-th9.1}, \ref{survey-th9.2} і~\ref{survey-th9.4}
встановлені в \cite{AnopMurach14UMJ7, Anop14Dop4}, а теорема~\ref{survey-th9.3}~--- в \cite[теорема~3]{AnopKasirenko16MFAT4}. Їх версії для різних класів еліптичних крайових задач доведені в
\cite{KasirenkoMurach18UMJ11} (задачі з крайовими умовами високих порядків), \cite{Anop13Collection2} (задачі в многозв'язній області) і \cite{Anop14Collection2} (матричні задачі), \cite[п.~3]{KasirenkoMurachChepurukhina19Dop3} (задачі на многовиді з краєм). В уточненій соболєвській шкалі такі класи задач досліджені в \cite{MikhailetsMurach06UMJ3, MikhailetsMurach07UMJ5, Murach07Dop4, Murach07Dop6} (відповідні результати викладено у монографії \cite[пп. 4.1, 5.3]{MikhailetsMurach14}). У праці Г.~Шлензак \cite{Slenzak74} доведено версію теореми~\ref{survey-th9.1} для більш вузької множини функціональних параметрів $\varphi$, яка не описана конструктивно. У~соболєвському випадку $\varphi(t)\equiv t^{s}$ теореми \ref{survey-th9.1}--\ref{survey-th9.3} добре відомі (принаймні, якщо $s\geq2q$, $\Omega_0:=\Omega$ і $\Gamma_{0}:=\Gamma$; див., наприклад, огляд \cite[\S\S~2, 4, 6]{Agranovich97}).

\section{Еліптичні крайові задачі з параметром}\label{survey-sect10}

Нехай цілі числа $q$ і $m_{1},\ldots,m_{q}$ такі як у розд.~\ref{survey-sect9}. В області $\Omega$ розглядається лінійна крайова задача
\begin{gather}\label{survey-4f1a}
A(\lambda)u(x)\equiv
\sum_{r=0}^{2q}\,\lambda^{2q-r}\sum_{|\mu|\leq r}a_{r}^{\mu}(x)D^{\mu}=f(x),\quad x\in\Omega,\\
B_{j}(\lambda)u(x)\equiv
\sum_{r=0}^{m_{j}}\,\lambda^{m_{j}-r}
\sum_{|\mu|\leq r}b_{j,r}^{\mu}(x)D^{\mu}=g_{j}(x),
\quad x\in\Gamma,\quad j=1,\ldots,q. \label{survey-4f1b}
\end{gather}
Тут $\lambda$~--- комплексний параметр, та усі коефіцієнти $a_{r}^{\mu}\in C^{\infty}(\overline{\Omega})$ і $b_{j,r}^{\mu}(x)\in C^{\infty}(\Gamma)$. Покладемо $B(\lambda):=(B_{1}(\lambda),\ldots,B_q(\lambda))$.

Пов'яжемо з диференціальними виразами $A(\lambda)=A(x,D,\lambda)$, де $x\in\overline{\Omega}$, і $B_{j}(\lambda)=B_{j}(x,D,\lambda)$, де $x\in\Gamma$, такі однорідні поліноми від $(\xi,\lambda)\in\mathbb{C}^{n+1}$:
\begin{equation*}
A^{\circ}(x,\xi,\lambda):=
\sum_{r=0}^{2q}\,\lambda^{2q-r}\sum_{|\mu|=r}a_{r}^{\mu}(x)\,\xi^{\mu}
\quad\mbox{і}\quad
B^{\circ}_{j}(x,\xi,\lambda):=\sum_{r=0}^{m_{j}}\,\lambda^{m_{j}-r}
\sum_{|\mu|=r}b_{j,r}^{\mu}(x)\,\xi^{\mu}.
\end{equation*}

Нехай $K$~--- деякий фіксований замкнений кут на комплексній площині з вершиною у початку координат (він може вироджуватися у промінь). Припускаємо, що крайова задача \eqref{survey-4f1a}, \eqref{survey-4f1b} є \textit{еліптичною з параметром} в куті $K$ в області $\Omega$, тобто задовольняє такі умови:
\begin{itemize}
\item[(i)] $A^{\circ}(x,\xi,\lambda)\neq0$ для кожної точки $x\in\overline{\Omega}$ та довільних $\xi\in\mathbb{R}^{n}$ і $\lambda\in K$ таких, що $|\xi|+|\lambda|\neq0$.
\item[(ii)] Для довільних точки $x\in\Gamma$, вектора $\tau\in\mathbb{R}^{n}$, дотичного до $\Gamma$ в точці $x$, і параметра $\lambda\in K$ таких, що $|\tau|+|\lambda|\neq0$, система многочленів $B^{\circ}_{j}(x,\tau+\zeta\nu(x),\lambda)$, $j=1,\ldots,q$, аргументу $\zeta\in\mathbb{C}$ є лінійно незалежною за модулем многочлена $\prod_{j=1}^{q}(\zeta-\zeta^{+}_{j}(x,\tau,\lambda))$; тут $\zeta^{+}_{j}(x,\tau,\lambda)$, де $j=1,\ldots,q$, є усі $\zeta$-корені многочлена $A^{\circ}(x,\tau+\zeta\nu(x),\lambda)$ з $\mathrm{Im}\,\zeta>\nobreak0$ (виписані з урахуванням їх кратності).
\end{itemize}

Умова (ii) коректна, оскільки з умови (i) випливає, що многочлен $A^{\circ}(x,\tau+\zeta\nu(x),\lambda)$ аргументу $\zeta\in\mathbb{C}$ має  $q$ коренів з $\mathrm{Im}\,\zeta>0$ і $q$ коренів з $\mathrm{Im}\,\zeta<0$ (з урахуванням їх кратності) \cite[твердження~2.2]{AgranovichVishik64}.

Оскільки, для кожного фіксованого $\lambda\in\mathbb{C}$ крайова задача \eqref{survey-4f1a}, \eqref{survey-4f1b} еліптична в області $\Omega$, то відображення $u\mapsto(A(\lambda)u,B(\lambda)u)$, де $u\in C^\infty(\overline{\Omega})$, продовжується єдиним чином до фредгольмового обмеженого оператора
\begin{equation}\label{survey-4f11}
(A(\lambda),B(\lambda)):H^{\varphi}(\Omega)\rightarrow
\mathcal{H}^{\varphi\varrho^{-2q}}(\Omega,\Gamma)
\end{equation}
для довільного $\varphi\in\mathrm{RO}$ такого, що $\sigma_{0}(\varphi)>m+1/2$.

\begin{theorem}\label{survey-th10.1}
Існує число  $\lambda_{1}>0$ таке, що для довільних $\lambda\in K$ з $|\lambda|\geq\lambda_{1}$ і $\varphi\in\mathrm{RO}$ з $\sigma_{0}(\varphi)>m+1/2$ оператор \eqref{survey-4f11} є ізоморфізмом. Більше того, для кожного $\varphi\in\mathrm{RO}$ з $\sigma_{0}(\varphi)>2q$ існує число $c=c(\varphi)\geq\nobreak1$ таке, що
\begin{gather*}
c^{-1}\bigl(\|u\|_{\varphi,\Omega}+
\varphi(|\lambda|)\|u\|_{\Omega}\bigr)\leq\\
\|f\|_{\varphi\varrho^{-2q},\Omega}+
\varphi(|\lambda|)|\lambda|^{-2q}\|f\|_{\Omega}+
\sum_{j=1}^{q}\bigl(\|g_{j}\|_{\varphi\varrho^{-m_j-1/2},\Gamma}+
\varphi(|\lambda|)|\lambda|^{-m_j-1/2}\|g_{j}\|_{\Gamma}\bigr)\leq\\
\leq c\,
\bigl(\|u\|_{\varphi,\Omega}+\varphi(|\lambda|)\|u\|_{\Omega}\bigr)
\end{gather*}
для кожного $\lambda\in K$ з $|\lambda|\geq\lambda_{1}$ і довільних
$u\in H^{\varphi}(\Omega)$ і $(f,g)\in\mathcal{H}^{\varphi\varrho^{-2q}}(\Omega,\Gamma)$, які задовольняють крайову задачу \eqref{survey-4f1a}, \eqref{survey-4f1b}.
Число $c$ не залежить від $\lambda$, $u$, $f$ і $g=(g_1,\ldots,g_q)$,
а $\|\cdot\|_{\Omega}$ і $\|\cdot\|_{\Gamma}$ позначають норми у просторах $L_2(\Omega)$ і $L_2(\Gamma,dx)$ відповідно.
\end{theorem}

Звідси випливає, що індекс оператора \eqref{survey-4f11} дорівнює нулю для будь-яких $\lambda\in\mathbb{C}$ і $\varphi\in\mathrm{RO}$ з $\sigma_{0}(\varphi)>m+1/2$.

Теорема~\ref{survey-th10.1} доведена в \cite[пп. 4,~6]{AnopMurach14MFAT2}, а для уточненої соболєвської шкали~--- в \cite[п.~7]{MikhailetsMurach07UMJ5} (див. також монографію \cite[п.~4.1.4]{MikhailetsMurach14}). Для гільбертових просторів Соболєва, залежних від параметра, вона встановлена М.~С.~Аграновичем і М.~І.~Вішиком \cite[\S\S~4,~5]{AgranovichVishik64} і застосована до доведення коректної розв'язності параболічних початково-крайових задач в анізотропних просторах Соболєва (див. також огляд \cite[п.~3.2]{Agranovich97}).

\section{Еліптичні задачі з грубими крайовими даними}\label{survey-sect11}

Тут розглядається регулярна еліптична крайова задача \eqref{survey-3f1},
\eqref{survey-3f2}, в якій праві частини $g_{1},\ldots,g_{q}$ крайових умов є довільними розподілами (як завгодно низької регулярності) на $\Gamma$, тобто, інакше кажучи, $g_{1},\ldots,g_{q}$ є грубими даними.
До такої задачі результати п.~\ref{survey-sect9} не застосовні, оскільки в них усі $g_{j}$~--- принаймні (щодо регулярності) квадратично інтегровні функції на $\Gamma$. Останню умову не задовольняють дельта-функції,  узагальнені похідні негладких функцій, функції зі степеневими особливостями, різні типи білого шуму на $\Gamma$; крім того, задачі з динамічними крайовими умовами приводять до задач з грубими даними (див., наприклад, \cite{DenkPlossRauSeiler23, FageotFallahUnser17, Hummel21BJMA, Roitberg96, Veraar11}).

Крайовий оператор $B_{j}$ не можна коректно означити на просторі $H^{(s)}(\Omega)$, якщо $s\leq m_{j}+1/2$. У такому випадку замість цього простору доводиться брати його частину, підпорядковану умові, яка полягає у тому, що $Au$ має вищу регулярність, ніж $Au\in H^{(s-2q)}(\Omega)$. У~соболєвському випадку, можна взяти умову $Au\in L_{2}(\Omega)$ як показали Ж.-Л.~Ліонс і Е.~Мадженес \cite{LionsMagenes62V, LionsMagenes63VI}, або дещо більш слабкі умови \cite{BehrndtHassideSnoo20, LionsMagenes72, MikhailetsMurach14, Murach09MFAT2}.

У рамках цього підходу уведемо для довільних параметрів $\alpha,\eta\in\mathrm{RO}$ гільбертів простір
\begin{equation*}
H^{\alpha}_{A,\eta}(\Omega):=
\bigl\{u\in H^{\alpha}(\Omega):\,Au\in H^\eta(\Omega)\bigr\},
\end{equation*}
наділений скалярним добутком і нормою
\begin{equation*}
(u_1,u_2)_{\alpha,A,\eta}:=
(u_1,u_2)_{\alpha,\Omega}+(Au_1,Au_2)_{\eta,\Omega},\quad
\|u\|_{\alpha,A,\eta}:=(u,u)_{\alpha,A,\eta}^{1/2};
\end{equation*}
тут $Au$ розуміємо у сенсі теорії розподілів. У соболєвському випадку,
коли функції $\alpha$ і $\eta$ степеневі цей простір досліджено в \cite{KasirenkoMikhailetsMurach19}.

Виберемо довільно функціональний параметр $\varphi\in\mathrm{RO}$, який задовольняє умову $\sigma_0(\varphi)\leq-1/2$. Крім того, виберемо дійсні числа $s_0$, $s_1$ і $\lambda$ такі, що $s_0<\sigma_0(\varphi)$, $s_1>\sigma_1(\varphi)$, $\lambda>-1/2$ і
\begin{equation*}
\left\{
  \begin{array}{ll}
    \lambda\leq s_{1},&\hbox{якщо}\;\;\;\sigma_1(\varphi)\geq-1/2;\\
    s_{1}<-1/2,&\hbox{якщо}\;\;\;\sigma_1(\varphi)<-1/2.
  \end{array}
\right.
\end{equation*}
Покладемо $\theta:=(s_1-\lambda)/(s_1-s_0)$ і означимо функцію $\eta(t)$  аргументу $t\geq1$ за формулою
\begin{equation*}
\eta(t):=\left\{
  \begin{array}{ll}
    t^{(1-\theta)s_{1}}\varphi(t^\theta),
    &\hbox{якщо}\;\;\;\sigma_1(\varphi)\geq-1/2;\\
    t^{\lambda},&\hbox{якщо}\;\;\;\sigma_1(\varphi)<-1/2.
  \end{array}
\right.
\end{equation*}
Маємо включення $\eta\in\mathrm{RO}$ і неперервне вкладення $H^{\eta}(\Omega)\hookrightarrow H^{\lambda}(\Omega)\cap H^{\varphi}(\Omega)$.

\begin{theorem}\label{survey-th11.1}
Множина $C^{\infty}(\overline{\Omega})$ щільна у просторі  $H^{\varphi\varrho^{2q}}_{A,\eta}(\Omega)$ і відображення $u\mapsto(Au,Bu)$, де $u\in C^\infty(\overline{\Omega})$, продовжується єдиним чином (за неперервністю) до обмеженого лінійного оператора
\begin{equation}\label{survey-th10.1-operator}
(A,B):H^{\varphi\varrho^{2q}}_{A,\eta}(\Omega)\to H^{\eta}(\Omega)\oplus
\bigoplus_{j=1}^{q}H^{\varphi\varrho^{2q-m_j-1/2}}(\Gamma)=:
\mathcal{H}^{\eta,\varphi}(\Omega,\Gamma).
\end{equation}
Цей оператор фредгольмів. Його ядро збігається з $N$, а область
значень складається з усіх векторів $(f,g_1,\ldots,g_q)\in\mathcal{H}^{\eta,\varphi}(\Omega,\Gamma)$, які задовольняють умову \eqref{survey-range}. Індекс оператора \eqref{survey-th10.1-operator} дорівнює $\dim N-\dim N^{+}$ й не залежить від~$\varphi$ і~$\eta$.
\end{theorem}

Якщо крайові дані $g_{1},\ldots,g_{q}$~--- довільно вибрані розподіли на $\Gamma$, то знайдеться число $s=s(g)\in\mathbb{R}$ таке, що
\begin{equation*}
\bigl(\varphi\in\mathrm{RO},\;\sigma_1(\varphi)\leq s\bigr)
\;\Longrightarrow\;
g=(g_{1},\ldots,g_{q})\in
\bigoplus_{j=1}^{q}H^{\varphi\varrho^{2q-m_j-1/2}}(\Gamma).
\end{equation*}
Отже, теорема~\ref{survey-th11.1} застосовна до цих даних, за умови, що права частина $f$ еліптичного рівняння має достатню регулярність у тому сенсі, що $f\in H^{\eta}(\Omega)$.

Як звичайно, $\mathcal{S}'(\Omega)$~--- лінійний топологічний простір звужень усіх розподілів $w\in\mathcal{S}'(\mathbb{R}^{n})$ в область $\Omega$. Якщо розподіл $u\in\mathcal{S}'(\Omega)$ задовольняє рівняння $Au=f$ для деякого $f\in H^{-1/2+}(\Omega)$, то $u$ належить області визначення оператора \eqref{survey-th10.1-operator}, де $\varphi(t)\equiv t^{s}$ і $\eta(t)\equiv t^{\lambda}$ для деяких чисел $s<-1/2$ і $\lambda>-1/2$. Тому вектор  $g:=Bu\in(\mathcal{D}'(\Gamma))^{q}$ коректно означений за замиканням на підставі теореми~\ref{survey-th10.1}, розглянутої у випадку просторів Соболєва. Отже, умови \eqref{survey-3f1},
\eqref{survey-3f2} мають сенс, якщо $u\in\mathcal{S}'(\Omega)$, $f\in H^{-1/2+}(\Omega)$ і $g\in(\mathcal{D}'(\Gamma))^{q}$. У цьому випадку розподіл $u$ називаємо \textit{узагальненим розв'язком} крайової задачі \eqref{survey-3f1}, \eqref{survey-3f2}.

Сформулюємо властивості цього розв'язку стосовно розширеної соболєвської шкали. Нехай множини $\Omega_{0}$ і $\Gamma_{0}$~--- такі самі як в розд.~\ref{survey-sect9}.

\begin{theorem}\label{survey-th11.2}
Припустимо, що розподіл $u\in\mathcal{S}'(\Omega)$ є узагальненим
розв'язком крайової задачі \eqref{survey-3f1}, \eqref{survey-3f2}, праві
частини якої задовольняють умови $f\in H^\eta_{\mathrm{loc}}(\Omega_{0},\Gamma_{0})\cap H^{-1/2+}(\Omega)$ і
$g_j\in H^{\varphi\varrho^{2q-m_j-1/2}}_{\mathrm{loc}}(\Gamma_0)$ для кожного номера $j\in\{1,\ldots,q\}$. Тоді $u\in
H^{\varphi\varrho^{2q}}_{\mathrm{loc}}(\Omega_0,\Gamma_0)$.
\end{theorem}

\begin{theorem}\label{survey-th11.3}
Припустимо, що розподіл $u\in\mathcal{S}'(\Omega)$ задовольняє умову теореми~$\ref{survey-th11.2}$. Довільно виберемо число $r>0$ і функції $\chi,\chi_{1}\in C^\infty(\overline{\Omega})$, підпорядковані умові \eqref{cutoff-functions}. Тоді існує число $c=c(\varphi,\eta,r,\chi,\chi_{1})>0$ таке, що
\begin{equation*}
\|\chi u\|_{\varphi\varrho^{2q},\Omega}\leq
c\biggl(\|\chi_{1}f\|_{\eta,\Omega}+
\sum_{j=1}^{q}\|\chi_{1}g_{j}\|_{\varphi\varrho^{2q-m_j-1/2},\Gamma}+
\|\chi_{1}u\|_{\varphi\varrho^{2q-r},\Omega}
\biggr),
\end{equation*}
причому $c$ не залежить від $u$ і $(f,g)=Au$.
\end{theorem}

\begin{theorem}\label{survey-th11.4}
Нехай ціле $\lambda\geq0$. Припустимо, що розподіл $u\in \mathcal{S}'(\Omega)$ задовольняє умову теореми~$\ref{survey-th11.2}$, причому
\begin{equation}\label{survey-th11.4-int-cond}
\int_1^{\infty}t^{2\lambda+n-1-4q}\varphi^{-2}(t)dt<\infty.
\end{equation}
Тоді $u\in C^{\lambda}(\Omega_{0}\cup\Gamma_{0})$.
\end{theorem}

Умова \eqref{survey-th11.4-int-cond} є точною у цій теоремі на класі зазначених розв'язків. Природно, для розподілу $u\in\mathcal{D}'(\Omega)$
включення $u\in C^{\lambda}(\Omega_{0}\cup\Gamma_{0})$ означає існування функції $u_0\in C^{\lambda}(\Omega_{0}\cup\Gamma_{0})$ такої, що
\begin{equation*}
\bigl(v\in C^{\infty}_{0}(\mathbb{R}^{n}),\;
\mathrm{supp}\,v\subset\Omega_{0}\bigr)\Longrightarrow
(u,v)_{\Omega}=\int_{\Omega_{0}}
u_0(x)\overline{v(x)}dx;
\end{equation*}
тут $(u,v)_{\Omega}$~--- значення розподілу $u$ на основній функції~$v$.

Теореми \ref{survey-th11.2}--\ref{survey-th11.4} зберігають силу для довільного параметра $\varphi\in\mathrm{RO}$ з $\sigma_0(\varphi)>-1/2$, якщо у них покласти $\eta:=\varphi$. Вони доведені в \cite[пп. 4,~6]{AnopDenkMurach21} (перша з них встановлена раніше в \cite[п.~4]{AnopMurach15Coll2} у випадку, коли $s_{1}\geq0$ і $\lambda=0$).

В уточненій соболєвській шкалі регулярні еліптичні крайові задачі з грубими даними досліджено в \cite{MikhailetsMurach11Dop4}, а для модифікації цієї шкали за Я.~А.~Ройтбергом~--- в \cite{MikhailetsMurach08UMJ4} (див. також монографію \cite[пп. 4.2, 4.5]{MikhailetsMurach14}). Для нерегулярних еліптичних крайових задач відповідні дослідження проведені в \cite{AnopKasirenkoMurach18UMJ3, KasirenkoMurach18MFAT2, Kasirenko18Dop2} (в останній праці використані показники гладкості класу $\mathrm{RO}$ з нульовими індексами Матушевської). У просторах Соболєва дослідження загальних еліптичних крайових задач з грубими даними започатковано Ж.-Л.~Ліонсом, E.~Мадженесом \cite{LionsMagenes62V, LionsMagenes63VI, LionsMagenes72} і Я.~А.~Ройтбергом \cite{Roitberg96, Roitberg65} (у цьому зв'язку дивитися \cite{MikhailetsMurach14, Murach09MFAT2}, недавні праці \cite{BehrndtHassideSnoo20, DenkPlossRauSeiler23, Hummel21JEE, KasirenkoMikhailetsMurach19} та наведену там бібліографію).

Звісно, теореми \ref{survey-th11.1}--\ref{survey-th11.4} застосовні
до задачі \eqref{survey-3f1}, \eqref{survey-3f2} у важливому випадку, коли еліптичне диференціальне рівняння \eqref{survey-3f1} однорідне, тобто розподіл $f=0$ в області~$\Omega$. Теорема про фредгольмовість відповідного оператора представляє окремий інтерес. Областю визначення цього оператора служить гільбертів простір \begin{equation*}
H^{\varphi}_{A}(\Omega):=
\bigl\{u\in H^{\varphi}(\Omega):Au=0\;\,\mbox{in}\;\,\Omega\bigr\},
\end{equation*}
наділений скалярним добутком з простору $H^{\varphi}(\Omega)$, де $\varphi\in\mathrm{RO}$. Оскільки диференціальний вираз $A$ еліптичний на $\overline{\Omega}$, то $H^{\varphi}_{A}(\Omega)\subset C^{\infty}(\Omega)$; проте $H^{\varphi}_{A}(\Omega)\not\subset C^{\infty}(\overline{\Omega})$. Позначимо через $C^{\infty}_{A}(\overline{\Omega})$ простір усіх функцій $u\in C^{\infty}(\overline{\Omega})$ таких, що $Au=0$ на $\overline{\Omega}$.

\begin{theorem}\label{survey-th11.5}
Для кожного $\varphi\in\mathrm{RO}$ множина $C^{\infty}_{A}(\overline{\Omega})$ щільна у просторі $H^{\varphi}_{A}(\Omega)$ і відображення $u\mapsto Bu$, де $u\in C^{\infty}_{A}(\overline{\Omega})$, продовжується єдиним чином (за неперервністю) до обмеженого лінійного оператора
\begin{equation}\label{survey-th11.1-f7.2}
B_{A}:H^{\varphi}_{A}(\Omega)\to
\bigoplus_{j=1}^{q}H^{\varphi\rho^{-m_j-1/2}}(\Gamma)=:
\mathcal{H}_{\varphi}(\Gamma).
\end{equation}
Цей оператор фредгольмів. Його ядро збігається з $N$, а область
значень складається з усіх векторів $(g_1,\ldots,g_q)\in\mathcal{H}_{\varphi}(\Gamma)$, які задовольняють умову
\begin{equation*}
\sum_{j=1}^{q}\,(g_{j},\,C^{+}_{j}v)_{\Gamma}=0
\quad\mbox{для довільного}\quad v\in N^{+}.
\end{equation*}
Індекс оператора \eqref{survey-th11.1-f7.2} дорівнює $\dim N-\dim N^{+}_1$, де $N^{+}_1$~--- скінченновимірний простір усіх векторів
$(C^{+}_{1}v,\ldots,C^{+}_{q}v)$ таких, що $v\in N^{+}$.
\end{theorem}

Звісно, $\dim N^{+}_{1}\leq\dim N^{+}$, причому можлива строга нерівність  \cite[теорема~13.6.15]{Hermander07v2}. Доповнимо цей результат достатньою умовою рівномірної збіжності послідовності розв'язків однорідного еліптичного рівняння та їх похідних.

\begin{theorem}\label{survey-th11.6}
Нехай ціле $\lambda\geq0$. Припустимо, що послідовність $(u_{k})_{k=1}^{\infty}\subset\mathcal{S}'(\Omega)$ задовольняє такі умови
\begin{itemize}
  \item[(i)] $Au_{k}=0$ в $\Omega$ для кожного $k\in\mathbb{N}$;
  \item[(ii)] ця послідовність збігається у соболєвському просторі  $H^{(-r)}(\Omega)$ для достатньо великого $r$;
  \item[(iii)] послідовність розподілів  $(B_{A}u_{k})\!\upharpoonright\!\Gamma_{0}$ збігається у просторі $\mathcal{H}_{\varphi}(\Gamma_{0})$ для деякого параметра $\varphi$, який задовольняє умову \eqref{survey-th9.4-cond}.
\end{itemize}
Тоді $u_{k}\in C^{\lambda}(\Omega\cup\Gamma_{0})$ для кожного $k\in\mathbb{N}$ і будь-яка послідовність $(D^{\mu}u_{k})_{k=1}^{\infty}$, де $|\mu|\leq\lambda$, рівномірно збігається на довільній замкненій (у просторі $\mathbb{R}^{n}$) підмножині множини $\Omega\cup\Gamma_{0}$.
\end{theorem}

Тут $\mathcal{H}_{\varphi}(\Gamma_{0})$~--- гільбертів простір звужень усіх вектор-функцій $g\in\mathcal{H}_{\varphi}(\Gamma)$ на відкриту підмножину $\Gamma_{0}$ межі $\Gamma$ (він означається аналогічно до просторів на $\Omega$).

Теореми \ref{survey-th11.5} і \ref{survey-th11.6} доведені в \cite[пп. 7.1, 7.2]{AnopDenkMurach21} (див. також \cite[пп. 4,~6]{AnopMurach16Coll2}). З~теореми~\ref{survey-th11.5} випливає \cite[п.~7.3]{AnopDenkMurach21}, що для простору $H^{\varphi}_{A}(\Omega)$ правильні аналоги інтерполяційних теорем $\ref{th-int-Sobolev}$--$\ref{th-int-space}$ (див. також \cite[п.~3]{AnopMurach18Dop3}). Для уточненої соболєвської шкали теорема~\ref{survey-th11.5} доведена в \cite[пп. 1,~6]{MikhailetsMurach06UMJ11} і  \cite[п.~3.3]{MikhailetsMurach14}.

Регулярні еліптичні крайові задачі з однорідними крайовими умовами і грубою правою частиною еліптичного рівняння досліджені в уточненій соболєвській шкалі в \cite{MikhailetsMurach05UMJ5, MikhailetsMurach06UMB4} і \cite[п.~3.4]{MikhailetsMurach14} (для соболєвських просторів відповідні результати викладено в монографіях \cite[розд.~III, \S~6, пп. 4--6]{Berezansky68} і \cite[\S~5.5]{Roitberg96}).

\section{Еліптичні задачі з додатковими невідомими функціями\\ у крайових умовах}\label{survey-sect12}

Виберемо довільно цілі числа $q\geq1$, $\varkappa\geq1$, $m_{1},\ldots,m_{q+\varkappa}$, і $r_{1},\ldots,r_{\varkappa}$. В області $\Omega$ розглядається лінійна крайова задача
\begin{gather}\label{survey-sect11-12f1}
Au=f\quad\mbox{в}\;\;\Omega,\\
B_{j}u+\sum_{k=1}^{\varkappa}C_{j,k}v_{k}=g_{j}\quad\mbox{на}\;\;\Gamma,
\quad j=1,...,q+\varkappa, \label{survey-sect11-12f2}
\end{gather}
де розподіл $u$ в $\Omega$ і $\varkappa$ розподілів $v_{1},\ldots,v_{\varkappa}$ на $\Gamma$ є шуканими. Тут
$A:=A(x,D)$~--- лінійний диференціальний оператор (ЛДО) на $\overline{\Omega}$, кожен $B_{j}:=B_{j}(x,D)$~--- крайовий ЛДО на $\Gamma$ і кожен $C_{j,k}:=C_{j,k}(x,D_{\tau})$~---  дотичний ЛДО на~$\Gamma$. Порядки цих операторів задовольняють умови $\mathrm{ord}\,A=2q$, $\mathrm{ord}\,B_{j}\leq m_{j}$ і $\mathrm{ord}\,C_{j,k}\leq m_{j}+r_{k}$, а їх коефіцієнти є функціями класу $C^{\infty}(\overline{\Omega})$ і $C^{\infty}(\Gamma)$ відповідно
(як звичайно ЛДО від'ємного порядку вважаємо рівним нулю.) Отже, порівняно з крайовою задачею \eqref{survey-3f1}, \eqref{survey-3f2}, задача \eqref{survey-sect11-12f1}, \eqref{survey-sect11-12f2} має $\varkappa$ додаткових невідомих функцій $v_{1},\ldots,v_{\varkappa}$ у крайових умовах. Такі задачі уперше розглянув Б.~Лаврук \cite{Lawruk63a} для еліптичних рівнянь.

Покладемо $m:=\max\{m_{1},\ldots,m_{q+\varkappa}\}$ і припустимо, що $m\geq-r_{k}$ для кожного $k\in\{1,\ldots,\varkappa\}$. Це припущення природне; справді, якщо $m+r_{k}<0$ для деякого $k$, то всі оператори  $C_{1,k}$,..., $C_{q+\varkappa,k}$ рівні нулю, тобто шуканий розподіл $v_{k}$ відсутній в крайових умовах \eqref{survey-sect11-12f2}. Допускаємо випадок $m\geq2q$.

Позначимо через $A^{\circ}(x,\xi)$, $B_{j}^{\circ}(x,\xi)$ і $C_{j,k}^{\circ}(x,\tau)$ головні символи ЛДО $A(x,D)$, $B_{j}(x,D)$ і $C_{j,k}(x,D_{\tau})$ відповідно; при цьому вважаємо, що останні два ЛДО мають формальні порядки $m_{j}$ і $m_{j}+r_{k}$. Отже, $A^{\circ}(x,\xi)$ і $B_{j}^{\circ}(x,\xi)$~--- однорідні поліноми аргументу $\xi\in\mathbb{C}^{n}$ порядку $2q$ і $m_{j}$ відповідно. Крім того, $C_{j,k}^{\circ}(x,\tau)$~--- однорідний поліном аргументу $\tau$ порядку $m_{j}+r_{k}$, де $\tau$~--- дотичний вектор до межі $\Gamma$ в точці~$x$.

Припускаємо, що крайова задача \eqref{survey-sect11-12f1}, \eqref{survey-sect11-12f2} є еліптичною в $\Omega$, тобто   \cite[означення~3.1.2]{KozlovMazyaRossmann97} оператор $A$ задовольняє умови (i), (ii), сформульовані в розд.~\ref{survey-sect9}, і ще таку умову: для кожної точки $x\in\Gamma$ і будь-якого ненульового вектора $\tau\in\mathbb{R}^{n}$, дотичного до $\Gamma$ у точці
$x\in\Gamma$, крайова задача
\begin{gather*}
A^{\circ}(x,\tau+\nu(x)D_{t})\theta(t)=0\quad\mbox{при}\;\;t>0,\\
B_{j}^{\circ}(x,\tau+\nu(x)D_{t})\theta(t)\big|_{t=0}+
\sum_{k=1}^{\varkappa}C_{j,k}^{\circ}(x,\tau)\lambda_{k}=0,\quad j=1,...,q+\varkappa,\\
\theta(t)\to0\quad\mbox{при}\quad t\rightarrow\infty
\end{gather*}
має лише нульовий розв'язок. Тут, функція  $\theta\in C^{\infty}([0,\infty))$ і числа $\lambda_{1},\ldots,\lambda_{\varkappa}\in\mathbb{C}$ є шуканими. Крім того, $A^{\circ}(x,\tau+\nu(x)D_{t})$ і $B_{j}^{\circ}(x,\tau+\nu(x)D_{t})$~--- диференціальні оператори відносно $D_{t}:=i\partial/\partial t$, які отримуємо, поклавши $\zeta:=D_{t}$ в поліномах $A^{\circ}(x,\tau+\nu(x)\zeta)$ і $B_{j}^{\circ}(x,\tau+\nu(x)\zeta)$ аргументу~$\zeta\in\mathbb{C}$.
Наведена умова еквівалента умові (iii) з розд.~\ref{survey-sect9}, якщо $\varkappa=0$.

Пов'яжемо із крайовою задачею \eqref{survey-sect11-12f1}, \eqref{survey-sect11-12f2} лінійне відображення
\begin{equation}\label{survey-sect11-12f3}
\begin{gathered}
\Lambda:(u,v_{1},...,v_{\varkappa})\mapsto
\biggl(Au,B_{1}u+\sum_{k=1}^{\varkappa}C_{1,k}\,v_{k},...,
B_{q+\varkappa}u+\sum_{k=1}^{\varkappa}C_{q+\varkappa,k}\,v_{k}\biggr)\\
\mbox{де}\quad u\in C^{\infty}(\overline{\Omega})\quad\mbox{і}\quad
v_{1},\ldots,v_{\varkappa}\in C^{\infty}(\Gamma),
\end{gathered}
\end{equation}
та гільбертові простори
\begin{gather*}
\mathcal{D}^{\varphi}(\Omega,\Gamma):=H^{\varphi}(\Omega)\oplus
\bigoplus_{k=1}^{\varkappa}H^{\varphi\varrho^{r_{k}-1/2}}(\Gamma),\\
\mathcal{E}^{\varphi\varrho^{-2q}}(\Omega,\Gamma):=
H^{\varphi\varrho^{-2q}}(\Omega)\oplus
\bigoplus_{j=1}^{q+\varkappa}H^{\varphi\varrho^{-m_{j}-1/2}}(\Gamma),
\end{gather*}
де $\varphi\in\mathrm{RO}$. Позначимо через $\|\cdot\|_{\varphi}'$ норму в $\mathcal{D}^{\varphi}(\Omega,\Gamma)$, а через $\|\cdot\|_{\varphi\varrho^{-2q}}''$~--- норму в $\mathcal{E}^{\varphi\varrho^{-2q}}(\Omega,\Gamma)$.

\begin{theorem}\label{survey-th12.1}
Нехай $\varphi\in \mathrm{RO}$ і $\sigma_{0}(\varphi)>m+1/2$. Тоді відображення \eqref{survey-sect11-12f3} продовжується єдиним чином (за неперервністю) до обмеженого лінійного оператора
\begin{equation}\label{survey-sect11-12f9}
\Lambda:\mathcal{D}^{\varphi}(\Omega,\Gamma)\rightarrow
\mathcal{E}^{\varphi\varrho^{-2q}}(\Omega,\Gamma).
\end{equation}
Цей оператор нетерів. Його ядро лежить у просторі $C^{\infty}(\overline{\Omega})\times(C^{\infty}(\Gamma))^{\varkappa}$ і разом з індексом не залежить від $\varphi$.
\end{theorem}

Область значень оператора допускає опис за допомогою формально спряженої задачі відносно спеціальної формули Гріна \cite[п.~3.1.3]{KozlovMazyaRossmann97}.

Позначимо через $\mathcal{D}^{m+1/2+}(\Omega,\Gamma)$ об'єднання усіх просторів $\mathcal{D}^{\varphi}(\Omega,\Gamma)$ таких, що $\varphi\in\mathrm{RO}$ і $\sigma_{0}(\varphi)>m+1/2$. Вектор
\begin{equation*}
(u,v):=(u,v_{1},\ldots,v_{\varkappa})\in\mathcal{D}^{m+1/2+}(\Omega,\Gamma)
\end{equation*} називаємо
(сильним) \textit{узагальненим розв'язком} крайової задачі \eqref{survey-sect11-12f1}, \eqref{survey-sect11-12f2} з правою частиною
\begin{equation*}
(f,g):=(f,g_{1},\ldots,g_{q+\varkappa})\in \mathcal{D}'(\Omega)\times(\mathcal{D}'(\Gamma))^{q+\varkappa},
\end{equation*}
якщо $\Lambda(u,v)=(f,g)$ для деякого оператора \eqref{survey-sect11-12f9}. Цей розв'язок має властивості подібні до властивостей розв'язків регулярної еліптичної крайової задачі. Наведемо їх; множини $\Omega_{0}$ і $\Gamma_{0}$~--- такі самі як в розд.~\ref{survey-sect9}.

\begin{theorem}\label{survey-th12.2}
Припустимо, що вектор $(u,v)\in\mathcal{D}^{m+1/2+}(\Omega,\Gamma)$ є узагальненим розв'язком крайової задачі \eqref{survey-sect11-12f1}, \eqref{survey-sect11-12f2}, праві частини якої задовольняють умови $f\in H^{\varphi\varrho^{-2q}}_{\mathrm{loc}}(\Omega_0,\Gamma_0)$ і $g_j\in H^{\varphi\varrho^{-m_j-1/2}}_{\mathrm{loc}}(\Gamma_0)$ для кожного номера $j\in\{1,\ldots,q+\varkappa\}$ і деякого параметра $\varphi\in\mathrm{RO}$ з $\sigma_0(\varphi)>m+1/2$. Тоді $u\in
H^{\varphi}_{\mathrm{loc}}(\Omega_0,\Gamma_0)$ і $v_{k}\in H^{\varphi\varrho^{r_k-1/2}}_{\mathrm{loc}}(\Gamma_{0})$ для кожного номера $k\in\{1,\ldots,\varkappa\}$.
\end{theorem}

\begin{theorem}\label{survey-th12.3}
Припустимо, що вектор $(u,v)\in\mathcal{D}^{m+1/2+}(\Omega,\Gamma)$
задовольняє умову теореми~$\ref{survey-th12.2}$. Довільно виберемо функції $\chi,\chi_{1}\in C^\infty(\overline{\Omega})$, підпорядковані умові \eqref{cutoff-functions}. Тоді існує число $c=c(\varphi,\chi,\chi_{1})>0$ таке, що
\begin{equation*}
\|\chi(u,v)\|_{\varphi}'\leq c\,\bigl(\,\|\chi(f,g)\|_{\varphi\varrho^{-2q}}''+
\|\chi_{1}(u,v)\|_{\varphi\varrho^{-1}}'\bigl),
\end{equation*}
причому $c$ не залежить від $(u,v)$ і $(f,g)=\Lambda(u,v)$.
\end{theorem}

\begin{theorem}\label{survey-th12.4}
Нехай ціле $\lambda\geq0$. Припустимо, що вектор $(u,v)\in\mathcal{D}^{m+1/2+}(\Omega,\Gamma)$ задовольняє умову теореми~$\ref{survey-th12.2}$, в якій функціональний параметр $\varphi$ задовольняє \eqref{survey-th9.4-cond}. Тоді
Тоді $u\in C^{\lambda}(\Omega_{0}\cup\Gamma_{0})$.
\end{theorem}

\begin{theorem}\label{survey-th12.5}
Нехай ціле $\lambda\geq0$, $k\in\{1,\ldots,\varkappa\}$ і $\Gamma_{0}\neq\emptyset$. Припустимо, що вектор $(u,v)\in\mathcal{D}^{m+1/2+}(\Omega,\Gamma)$ задовольняє умову теореми~$\ref{survey-th12.2}$, в якій
\begin{equation*}
\int_1^{\infty}t^{2(\lambda-r_{k})+n-1}\varphi^{-2}(t)dt<\infty.
\end{equation*}
Тоді $v_{k}\in C^{l}(\Gamma_{0})$.
\end{theorem}

Ці теореми доведені в \cite{AnopChepurukhinaMurach21Axioms}. Їх окремі випадки встановлені в \cite{MurachChepurukhina20Dop8, Chepurukhina15Coll2} для розширеної соболєвської шкали і в \cite{KasirenkoChepurukhina17Coll2, Chepurukhina14Coll2} для уточненої соболєвської шкали. Задачі з грубими даними досліджено в \cite{ChepurukhinaMurach15MFAT1, ChepurukhinaMurach20MFAT2, MurachChepurukhina15UMJ5} в уточненій соболєвській шкалі і в \cite{ChepurukhinaMurach15MFAT1} для модифікації за Ройтбергом цієї шкали, а також в \cite{Anop19Dop2, Chepurukhina15Dop7} для однорідних еліптичних диференціальних рівнянь. У модифікованих за Ройтбергом просторах Соболєва такі задачі досліджено в монографіях \cite[розд. 3,~4]{KozlovMazyaRossmann97} і \cite[розд.~2]{Roitberg99}, а
у немодифікованих соболєвських просторах~--- в \cite{Chepurukhina16Coll2}.

\medskip

\textit{Ця робота фінансована Національною академією наук України. Автори були також підтримані грантом імені Марії Склодовської-Кюрі №~873071 за програмою Європейського Союзу Горизонт 2020~--- Рамкова програма з досліджень та інновацій: Спектральна оптимізація: від математики до фізики і передових технологій. (The authors were also supported by the European Union's Horizon 2020 research and innovation programme under the Marie Sk{\l}odowska-Curie grant agreement No~873071 (SOMPATY: Spectral Optimization: From Mathematics to Physics and Advanced Technology).)}

\end{document}